\documentclass[11pt]{amsart}
\usepackage{latexsym,amssymb,amsmath}
\usepackage{youngtab}
\input epsf
\DeclareMathOperator{\GL}{GL}

\newtheorem{theorem}{Theorem}[section]
\newtheorem{definition}[theorem]{Definition}

\newtheorem{proposition}[theorem]{Proposition}
\newtheorem{lemma}[theorem]{Lemma}
\newtheorem{corollary}[theorem]{Corollary}

\newcommand{\g}{\mathfrak{g}}

\newcommand{\ft}{\mathfrak{t}}

\newcommand{\fsl}{\mathfrak{sl}}
\newcommand{\fsp}{\mathfrak{sp}}
\newcommand{\fgl}{\mathfrak{gl}}
\newcommand{\fso}{\mathfrak{so}}

\def\tdeg{\text{deg} }

\def\tdim{\text{dim}\,}

\def\BC{\mathbb C}\def\BO{\mathbb O}\def\BS{\mathbb S}
\def\BA{\mathbb A}\def\BR{\mathbb R}\def\BH{\mathbb H}
\def\BP{\mathbb P}
\def\pp#1{\mathbb P^{#1}}

\def\fgl{\mathfrak g\mathfrak l}
\def\ta{\tilde{\alpha}}

\def\uG{\overline G}

\def\trace{{\rm trace}}

\def\cO{{\mathcal O}}
\def\cZ{{\mathcal Z}}
\def\cH{{\mathcal H}}
\def\CC{\mathbb C}
\def\RR{\mathbb R}
\def\HH{\mathbb H}
\def\AA{{\mathbb A}}
\def\BB{{\mathbb B}}
\def\OO{\mathbb O}

\def\ZZ{\mathbb Z}
\def\SS{\mathbb S}

\def\11{\mathbf 1}
\def\PP{\mathbb P}
\def\QQ{\mathbb Q}
\def\FF{\mathbb F}

\def\fh{{\mathfrak h}}

\def\fc{{\mathfrak c}}

\def\fsl{{\mathfrak {sl}}}
\def\fsp{{\mathfrak {sp}}}

\def\fso{{\mathfrak {so}}}
\def\fe{{\mathfrak e}}
\def\ff{{\mathfrak f}}

\def\fg{{\mathfrak g}}

\def\ft{{\mathfrak t}}

\def\l{\lambda}
\def\a{\alpha}
\def\ta{\tilde{\alpha}}
\def\o{\omega}

\def\b{\beta}
\def\g{\gamma}
\def\s{\sigma}

\def\th{\theta}

\def\up#1{{}^{({#1})}}
\def\e{\varepsilon}
\def\ot{{\mathord{\,\otimes }\,}}
\def\op{{\mathord{\,\oplus }\,}}

\def\lra{{\mathord{\;\longrightarrow\;}}}
\def\ra{{\mathord{\;\rightarrow\;}}}

\def\we{{\mathord{{\scriptstyle \wedge}}}}

\def\tr{{\rm trace}\;}
\def\dim{{\rm dim}\;}
\def\La{\Lambda}

\def\trace{{\rm trace}}
\def\cA{{\mathcal A}}
\def\cJ{{\mathcal J}}

\def\cO{{\mathcal O}}
\def\CC{\mathbb C}
\def\RR{\mathbb R}
\def\HH{\mathbb H}
\def\AA{{\mathbb A}}
\def\BB{{\mathbb B}}
\def\OO{\mathbb O}

\def\ZZ{\mathbb Z}
\def\SS{\mathbb S}

\def\11{\mathbf 1}
\def\PP{\mathbb P}
\def\QQ{\mathbb Q}
\def\FF{\mathbb F}

\def\fh{{\mathfrak h}}

\def\fsl{{\mathfrak {sl}}}
\def\fsp{{\mathfrak {sp}}}

\def\fso{{\mathfrak {so}}}
\def\fe{{\mathfrak e}}

\def\ff{{\mathfrak f}}

\def\fg{{\mathfrak g}}

\def\ft{{\mathfrak t}}

\def\l{\lambda}
\def\a{\alpha}
\def\ta{\tilde{\alpha}}
\def\o{\omega}

\def\b{\beta}
\def\g{\gamma}
\def\s{\sigma}

\def\th{\theta}

\def\up#1{{}^{({#1})}}
\def\e{\varepsilon}
\def\ot{{\mathord{\,\otimes }\,}}
\def\op{{\mathord{\,\oplus }\,}}

\def\lra{{\mathord{\;\longrightarrow\;}}}
\def\ra{{\mathord{\;\rightarrow\;}}}

\def\we{{\mathord{{\scriptstyle \wedge}}}}

\def\tr{{\rm trace}\;}
\def\dim{{\rm dim}\;}
\def\La#1{\Lambda^{#1}}

\def\xadg{X^{ad}_G}
\def\uG{\overline G}

\def\ug{\overline \fg}

\newcommand\rem{{\medskip\noindent {\em Remark}.}\hspace{2mm}}

\begin{document}
\title{The sextonions and $E_{7\frac 12}$ }
\author{J. M. Landsberg${}^1$, L. Manivel}
\footnote{Supported by NSF grant DMS-0305829}
\date{February 2004}

  \begin{abstract} We fill in the \lq\lq hole\rq\rq\ in the exceptional series
  of Lie algebras that was observed by Cvitanovic, Deligne, Cohen and deMan.
More precisely,
we show that the intermediate Lie algebra between $\fe_7$ and $\fe_8$  satisfies
some of the decomposition and dimension formulas of the exceptional simple Lie
algebras. A key role is played by the {\it sextonions}, a six dimensional
algebra between the quaternions and octonions. Using the sextonions, we
show simliar results
hold for the rows of an   expanded Freudenthal magic chart. We also obtain
new interpretations of the adjoint variety of the exceptional group $G_2$.
  \end{abstract}

\maketitle

\section{Introduction}

In \cite{del1,cdm, LMadv}
remarkable dimension formulas for the exceptional series of complex
simple Lie algebras were established, parametrizing the series by the
dual Coexeter number in \cite{del1,cdm} and using
the dimensions of composition algebras in
\cite{LMadv}. Cohen and deMan observed that all parameter values giving
rise to integer outputs in all the formulas of \cite{del1,cdm} were already
accounted for with  essentially one exception, which, were it the dimension
of a composition algebra, would be of dimension six and sit between the quaternions
and octonions to produce a Lie algebra sitting between $\fe_7$ and $\fe_8$.
B. Westbury brought this to our attention and pointed out that were this
the case, one would gain an entire new row of Freudenthal's magic chart.
We later learned
that this algebra,    which we call the {\it sextonion algebra}, 
had been observed earlier
as a curiosity \cite{jeur, klei}.

 In this paper we discuss the sextonions
and the extra row of the magic chart it gives rise to.  Along the way,
we discuss intermediate Lie algebras in general and their homogeneous
varieties, in particular the exceptional Lie algebra 
$\fe_{7\frac 12}$ defined by 
the triality construction of \cite{LMadv} applied to the sextonions. 
This Lie algebra is intermediate between $\fe_7$ and $\fe_8$. Of course
it is not simple but remarkably, shares most of the properties of the simple 
exceptional Lie algebras discovered by Vogel and Deligne. More generally,    
many of the dimension fomulas 
of \cite{del1,cdm, LMadv} are satisfied by the intermediate 
Lie algebras and
some of the decomposition formulas of \cite{LMseries} hold as well.

Let $\ug$ be a complex simple Lie algebra equipped with its
adjoint ($5$-step) grading induced by the highest root $\ta$:
$$
\ug = \ug_{-2}\op \ug_{-1}\op \ug_0 \op \ug_1\op \ug_2.
$$
Here $\ug_{ 2}\simeq \BC$ is the root space of $\ta$ and $\ug_0$ is reductive with
a one-dimensional center (except in type $A$ where the center is
two dimensional). Let $\fh = [\ug_0,\ug_0]$ be its semi-simple part. 

Introduce the {\it intermediate Lie algebra}
$$
\fg = \fh \op \ug_1\op \ug_2.
$$

Intermediate Lie algebras (sometimes in the  forms
$\fg '= \ug_0 \op \ug_1\op \ug_2$,
$\fg ''= \ug_0 \op \ug_1 $, $\fg '''= \fh \op \ug_1 $)
have appeared in \cite{shtepin, proc, GZ, GZ2}.
Shtepin used them to help decompose $\ug$-modules as $\fh$-modules
in a multiplicity free way to make Gelfand-Tsetlin bases.
Gelfand and Zelevinsky used them to make representation models
for the classical groups (and we hope the varieties discussed
here might lead to similar models for the exceptional groups, or
even that the triality model will give rise to uniform representation
models for all simple Lie groups).
Proctor made a detailed study of certain representations of the {\it odd 
symplectic Lie algebras} and proved a Weyl dimension formula for these.

% \medskip
 
\subsection*{Overview} In \S 2 we define
 and discuss the {\it adjoint varieties} of the intermediate Lie algebras.
 In \S 3 we give geometric interpretations of the adjoint variety of the
 exceptional group $G_2$, in particular
 we show that  it parametrizes sextonionic
 subalgebras of the octonions. We also give a new description of the
 variety of quaternionic subalgebras of the  octonions.
 In \S 4 we review the triality construction of Freudenthal's magic square
 and show how it applies to the sextonions. In \S 5, we discuss highest
 weight modules of intermediate Lie algebras, showing how to decompose
 the Cartan powers of the adjoint representation as an $\fh$ module. We also
 remark that some of Vogel's universal decomposition fomulas hold. In \S 6 
 we show that some of the more refined decomposition fomulas of
 \cite{LMseries} hold for the rows of the extended magic square. 
 Corresponding dimension formulas are stated and proven in \S 7.
Finally in \S\ref{geomsect} we describe the geometry
of closures of the orbits of highest weight
vectors  $\BP\overline{ (G.v) }\subset \BP V$
inside the preferred representations described in \S\ref{rowbyrowsect}.
In particular we get a new (slightly singular) Severi variety that 
we study   in detail. 

\medskip

\noindent{\it Notation}: We use the ordering of roots as in \cite{bou}.
Unless otherwise specified, all groups $G$ associated to a Lie
algebra $\fg$ are the adjoint groups.

\medskip

\noindent{\it Acknowledgement}: We would like to thank B. Westbury for
useful conversations and sharing his preprint \cite{BWW} with us.

\section{Adjoint varieties}

Let $\fh$ be a complex simple Lie algebra. The {\it adjoint variety}
$X^{ad}_H\subset \BP\fh$ is
the closed $H$-orbit in $\PP\fh$, where $H$ denotes the adjoint group of $\fh$.
 The adjoint variety parametrizes the highest root spaces:
given a line $\ell$ in $\fh$  which corresponds to a point of $X^{ad}_H$,
we can chose a Cartan subalgebra of $\fh$ and a set of positive roots, such that
$\ell=\fh_{\ta}$, the root space of the highest root $\ta$.
The goal of this section is to define a cousin of the adjoint variety
for intermediate Lie algebras, which we will also call the adjoint variety.

Up to the center, the intermediate Lie algebra $\fg$ coincides with the parabolic
subalgebra of $\ug$ which stabilizes the line $\ug_{\ta}\in X^{ad}_{\overline G}$.
Its reductive part  $\fh=[\ug_{0},\ug_{0}]$ is   simple exactly when the support of the adjoint
representation is an end of the Dynkin diagram of $\fg$. When $\fg\simeq\fsl_{n+1}$
is of type $A$, then $\fh=\fgl_{n-1}$ and we define the adjoint variety as
the closed $PGL_{n-1}$-orbit $\FF(1,n-2)\subset\PP\fgl_{n-1}$.
When $\fg\simeq\fso_m$ is of type $B$ or $D$, then $\fh\simeq\fsl_2\times\fso_{m-4}$
is not simple and it is not clear how to define the adjoint variety of $\fh$.
In fact we can take either the adjoint variety of $\fsl_2$, a plane conic,
or the adjoint variety $G_Q(2,m-4)$ of $\fso_{m-4}$. Notwithstanding this
difficulty for the $\fso_m$-case, we make the following:

\begin{definition} Let $\fg$ be an intermediate Lie algebra. The adjoint variety
$X^{ad}_G \subset \PP\fg $ of $\fg$ is the
closure   of the $G$-orbit of a highest weight line of  $\fh$.
\end{definition}

Recall that the $\fh$-module $\ug_1$ has very nice properties. It is 
{\it simple} (except in type $A$, for which we have two simple modules exchanged 
by an outer automorphism), and {\it minuscule}. 
The Lie bracket on $\fg$ induces an invariant symplectic form $\o\in\La 2\ug_1^*$ 
(canonically defined only up to scale), so in particular
$\fh\subset S^2\ug_1^*$. In fact $\fh$
generates the ideal of the closed orbit $H/Q\subset \BP \ug_1$
 which is Legendrian (in type $A$  it is the
union of two disjoint linear spaces, each of which is
Legendrian). See, e.g.,   \cite{LM1} for proofs of the above assertions.

Given $x\in \fh$, let $q^x\subset S^2\ug_1^*$ denote the
quadratic form it determines,
defined by $$q^x(v,w)=\frac 12\o (x.v,w)=\frac 12\o (x.w,v).$$
The  linear span $\langle q^x\rangle\subset \ug_1$ of this quadratic form 
is the image of the endomorphism  $L_x$ of $\ug_1$ given by the action of $x$.
Since $q^x(v,w)$ only depends on $x.v$ and $x.w$, we have an induced quadratic form 
on $\langle q^x\rangle$ defined by $q_x(v)=q^x(u,u)$ when $v=xu$.  

\begin{proposition} We have the following description of $\xadg$:
$$X^{ad}_G
=\overline{\{ [x, v,q_x(v) ]\mid
x\in \hat X^{ad}_H,\ v\in \langle \hat q^x\rangle\} }
\subset \BP \fg =\BP (\fh \op \ug_1\op \ug_2).
$$
Its affine tangent spaces at $p=[x,0,0]\in X^{ad}_H$ and $p_0=[0,0,1]$ are 
$$\hat T_pX^{ad}_G=\hat T_xX^{ad}_H \op\langle q^x\rangle, \qquad 
\hat T_{p_0}X^{ad}_G= \ug_1 \op \ug_2.$$
\end{proposition}

\proof Given $x\in\hat X^{ad}_H\subset \hat X^{ad}_G$, consider the action  of
$\exp(-u)$ on   $x$, for $u\in\ug_1$:
$$\exp(-u)x=x+[x,u]+\frac12 [u,[u,x]]=x+x.u+q^x(u).$$
The first claim follows, and 
the description of the tangent space at $p=[x,0,0]$ is clear.

Since $p_0$ is killed by $\fg$, $\hat T_{p_0}X^{ad}_G$ must be a $\fg$-submodule
of $\ug$,  contained in $\ug_1\op\ug_2$ since a linear action cannot move
a vector two steps in a grading. 
Since $\ug_1$ is irreducible (including in type $A$ if we take into account the 
$\ZZ_2$-action), and since the affine tangent space cannot be reduced to $\ug_2$,
there must be equality.
\qed

\medskip
Recall that the complex simple Lie algebras can be parametrized by their {\it Vogel 
parameters} $\a, \b, \g$ (roughly the Casimir eigenvalues of the nontrivial components of their
symmetric square, see \cite{vog2, LMuniv}). We normalize $\a=-2$, so that 
$t=\a+\b+\g=\check h$ is the dual Coxeter number. We distinguish
$\b$ from $\g$ as in \cite{LMuniv}. 

\begin{lemma}
For $x\in X^{ad}_H$, the dimension of $\langle q^x\rangle$ is equal to $\b$. 
\end{lemma}

\proof Let $x$ be the root space $\fg_{\ta}$ defined by a maximal root 
$\ta$ of $\fh$. Then the image of $L_x$ is the direct sum of the root spaces 
$\fg_{\mu+\ta}$, for $\mu$ a root of $\ug_1$ such that $\mu+\ta$ is again a root. 
By \cite{LMuniv}, Corollary 5.2, there exists exactly $\b$ such roots. \qed 

\begin{corollary}
The dimensions of the adjoint varieties of $\fh$,$\fg$ and $\ug$ are related
as follows:
$$
\begin{array}{lcl}
\dim  X^{ad}_{H}  & = &  2 \check h -3-2\b , \\
\dim  X^{ad}_{G}  & = &  2\check h -3- \b, \\
\dim  X^{ad}_{\overline G}  &  = &  2\check h -3
\end{array}
$$
unless $\ug=\fg_2$, in which case the formula holds with $2$ instead of $\b$.
\end{corollary}

\proof
The third equality was observed in \cite{knop}.
That the first and second lines differ by $\b$ follows
immediately from the fact that $\tdim \langle q^x\rangle =\b$,
and the description we gave of the tangent space to $X_G^{ad}$ at 
$[x,0,0]$, which is a generic point. 
Finally, the additional dimensions of the tangent space to
$X^{ad}_{\uG}$ at that point arise from the action of $\ug_{-1}$,
which is symmetric with the action of $\ug_1$ and thus contributes 
the same value $\b$. \qed

%\begin{remark}
%\rem  Similar to the $\b-2$ dimensional quadrics
%on $X^{ad}_{\ug}$ discussed above, the variety
%  $X_Q\subset \BP \ug_Q=\BP Y_2'$ parametrizes
%a family of $\b$-dimensional
%quadric sections of $X^{ad}_{\uG}$, that is, linear sections
 % $L\cap X^{ad}_{\uG}$ isomorphic to a smooth $\b$-dimensional quadric,
 % see \cite{LMseries}
%\end{remark}

\begin{corollary}
The adjoint  variety of an intermediate Lie algebra $\xadg$ is smooth if and only if
$\b (\ug) =1$. In general, the maximal excess dimension of its Zariski tangent spaces
is $\b (\ug) -1$.
\end{corollary}

\proof By semi-continuity, the excess dimension of the Zariski tangent space 
must be maximal at the point $p_0$ of $X^{ad}_G$. We have calculated that 
the dimension of the tangent space at that point is $2 \check h- 2$ while 
$\dim  X^{ad}_{G}= 2\check h -3- \b$.  \qed

\medskip\noindent {\it Example 1}.
 If $\ug=\fsp_{2n+2}$, the adjoint variety is $v_2(\PP^{2n+1})$. If $\ell$ is a point
of this variety, i.e., a line in $\CC^{2n+2}$, the Lie algebra $\fh$
may be identified with
$\fsp (V)$, for $V\simeq  \ell^{\perp}/\ell $, a vector space of dimension $2n$ 
endowed with the restriction of the original symplectic form, which is again 
symplectic. Its adjoint variety is $v_2(\PP V)$. The intermediate Lie algebra
$$\fg=\fsp(V)\op V\op \CC$$
is an {\it odd symplectic Lie algebra} (see \cite{proc}), and the corresponding adjoint 
variety is $v_2(\PP\ell^{\perp})$. We thus get a smooth variety with only two $G$-orbits,
the point $\ell$ and its complement.

\medskip\noindent {\it Example 2}.
If $\ug=\fsl_{n+1}$, the adjoint variety is $\FF_{1,n}$. A point of  this
variety is a pair $(\ell_0,H_0)$, with $\ell_0$ a line, and $H_0$ a hyperplane
containing $\ell_0$. The Lie algebra $\fh$ may be
identified  with $\fgl (V)$, where $V=H_0/\ell_0$ once we have
chosen a decomposition of $\CC^{n+1}$ as $\ell_0\op V\op\ell_1$, with $\ell_0\op V=H_0$.
Its adjoint variety is the set of pairs $(\ell\subset H)$, with $\ell$ a line
and $H$ a hyperplane  in $V$, defined by a linear form that we extend by zero on
$\ell_0\op\ell_1$. The intermediate Lie algebra is $$\fg=\fsl(V)\op V\op V^*\op\CC,$$ whose 
adjoint variety is 
$$\{(\ell,H)\in\FF_{1,n}, \quad \ell\subset H_0, \; \ell_0\subset H\}.$$
This variety has four $G$-orbits, and a unique singular point $(\ell_0,H_0)$, which is
a simple quadratic singularity.

\medskip
Let $\th : S^2\ug_1\ra \fh$ denote the projection map which is dual to the natural
inclusion $\fh\subset S^2\ug_1^*$. The duality on $\fh$ is taken with respect to the 
restriction of the Killing form $K$ of $\ug$. Explicitly,
\begin{eqnarray}
K(\theta(u),y)=\o(yu,u) \qquad {\rm for}\; u\in\ug_1, \; y\in\fh.
\end{eqnarray}

\begin{proposition}
Normalize the Killing form so that $(\a_0,\a_0)=2$.\hfill\hfill\linebreak
Let $x\in \hat X^{ad}_H$, $y\in\fg_0$ and $u\in\ug_1$. Then
\begin{eqnarray}
 xyxu &= &K(x,y)xu, \\
\theta (xu) &= &q^x(u,u)x.
\end{eqnarray}
\end{proposition}

\proof By homogeneity, we may suppose that $x=X_{\ta}$ belongs to the root space
$\fg_{\ta}$. The identity (1) being linear in
$y$ and $u$, we can let $u=X_{\theta}$ for some root $\theta\in\Phi_1$. If $y$ belongs to
the Cartan subalgebra, then $xu$ is an eigenvector of  $y$, thus $xyxu$ is a multiple of
$x^2u$, hence zero. Since $K(x,y)$ is also zero, we are done.

Now suppose that $y=X_{\s}$
is a root vector in $\fh$. Then $K(x,y)\ne 0$ if and only if $\s=-\ta$.
Recall that $\ug_1$ is a minuscule $\fh$-module, so that a root of $\ug_1$ is 
of the form $\g=\o_0+\chi$ with $-1\le \chi(H_{\tau})\le 1$ for every root $\tau$ of $\fh$. 
Moreover, for $u=X_{\g}$,
$X_{\tau}u\ne 0$ implies that $\chi (H_{\tau})=-1$. This implies that if $xyxu$ is nonzero,
$$-1=\chi (H_{\ta})=(\chi+\ta)(H_{\s})=(\chi+\ta+\s)(H_{\ta}),$$
hence $\s(H_{\ta})=-\ta(H_{\ta})=-2$. But this is possible only if $\s=-\ta$, in which case
$[y,x]=tH_{\ta}$ for some scalar $t\ne 0$, and
$$xyxu=x[y,x]u=t\g(H_{\ta})xu=t\chi(H_{\ta})xu=-txu.$$
With our normalization, $2t=tK(H_{\ta},H_{\ta})=K(H_{\ta},[y,x])=K([H_{\ta},y],x)=-2K(y,x)$,
thus $K(x,y)=-t$ and finally, $xyxu=K(x,y)u$, which is what we wanted to prove.

\smallskip
The second identity is an immediate consequence: from the equation defining $\theta$,
we get
$$K(y,\theta(v))=\o(v,yv)\qquad \forall y\in\fh, \,\forall v\in\ug_1.$$
For $x\in \hat X^{ad}_H$, we get using (1),
$$\begin{array}{rcccl}
K(y,\theta(xu)) & = & \o(yxu,xu)
 & = & -\o(xyxu,u) \\
 & = & -K(y,x)\o(xu,u)
 & = & -K(y,q^x(u,u)x),
\end{array}$$
as claimed. The proof is complete. \qed

\medskip\noindent Recall that the cone over the 
closed $H$-orbit in $\PP\ug_1$ is the set of vectors 
$v\in\ug_1$ such that 
$\theta(v)=0$, a space of quadratic equations parametrized by $\fh$ (see \cite{LM1}).

The following fact was observed case by case in \cite{LMuniv}:

\begin{corollary}
The adjoint variety parametrizes a family of $(\b-2)$-dimensional quadrics on
the closed $H$-orbit in $\PP\ug_1$.
\end{corollary}

\proof For each $x\in X^{ad}_H$, the linear space $\langle q^x\rangle\subset\PP\ug_1$
has dimension $\b$, and the identity $\theta(xu)=q^x(u)x=q_x(xu)x$ shows that on 
$\langle q^x\rangle$, the condition $\theta(v)=0$
reduces to a single quadratic condition. This means that the closed orbit cuts $\langle q^x\rangle$
along a quadric hypersurface.\qed

\medskip
In the case where $\fg=\fso_m$, we have two families of maximal quadrics in
the adjoint variety $G_Q(2,m)$, of dimension $4$ and $m-4$. We can either choose
$\b=4$ and $\g=m-4$, which corresponds to the component $G_Q(2,m-4)$ of $X_H^{ad}$,
or   $\b=m-4$ and $\g=4$, which corresponds to the one-dimensional
component $v_2(\PP^1)$.

\section{The sextonions}

\def\xadg2{X^{ad}_{G_2}}
Consider the adjoint variety $X^{ad}_{G_2}\subset \BP\fg_2$.
Since $G_2\subset SO(7, \underline q)$, we have $\xadg2\subset 
G_{\underline q}(2,Im\BO)$, the Grassmannian
of $\underline q$-isotropic planes in $\BC^7=Im\BO$. So the adjoint variety
is a $G_2$-invariant set of isotropic planes. Here
$\underline q$ is the restriction of the quadratic form 
$q(a,b)=Re(a\overline b)$
on $\BO$ to $Im\BO$, where  $Re(x):=\frac 12(x+\overline x)$
in $\BO$. 

We will say a plane $E\in G(2,Im\BO)$ is {\it null} if for all
$u,v\in E$, $uv=0$.

\begin{theorem}
The adjoint variety $\xadg2$, the closed $G_2$-orbit in $\PP\fg_2$, parametrizes:
\begin{enumerate}
\item null-planes in $\OO$;
\item rank two derivations of $\OO$, up to scalars;
\item six-dimensional subalgebras of $\OO$.
\end{enumerate}
\end{theorem}

\proof The correspondence between these three   objects is as follows: if $U$
is a null-plane in $\OO$, the orthogonal space
$V$ is a six-dimensional subalgebra, and there
is, up to scale, a unique skew-symmetric endomorphism of $\OO$ whose image is $U$
and kernel is $V$.

\smallskip
We first claim that $\xadg2= G(2,Im\OO)\cap\PP\fg_2$.
Note that this intersection is highly non-transverse
(of codimension $5$ in $G(2,7)$ instead of the expected $7$), although the set theoretic
intersection is a smooth variety. We   therefore use a
direct geometric
description in terms of the {\it associative form}  $\phi\in \La 3 Im\OO^*$ defined by
$$\phi(x,y,z)=Re [(xy)z-(zy)x].$$

  Bryant showed that the stablilizer of $\phi$ is exactly the group $G_2$,
see \cite{harvey}. Note that since
$\dim \La {3} Im\OO=\dim\fgl_7-\dim\fg_2$, the $GL_7$-orbit of $\phi$ in $\La 3 Im\OO$
is a dense open subset.

On the Grassmannian $G(2,Im\OO)$, we have a tautological rank two vector bundle $T$,
and a quotient bundle $Q$ of rank $5$. Consider the homogeneous vector bundle
 $E=Q^*\ot\La 2T^*$, of
rank $5$. By the Borel-Weil theorem, the space of global sections
of this vector bundle is $\Gamma (G(2,Im\OO),E)=\La 3(Im\OO)^*$. We can therefore
interpret $\phi$ as a generic section $\s$ of the vector bundle $E$, which is globally
generated, being irreducible as a homogeneous vector bundle. By Bertini, the zero-locus
of $\s$ is, if not empty, a smooth codimension $5$ subvariety of $G(2,Im\OO)$, hence
a $5$-dimensional smooth variety, $G_2$-invariant since $\phi$ is $G_2$-invariant.
But the adjoint variety $\xadg2$ is the $G_2$-orbit of minimal dimension, and this
dimension is five. So $\xadg2$ must be equal to the zero-locus of $\s$.

What is this zero-locus explicitly? If we choose a basis $u_1, u_2$ of a plane
$U$ in $Im\OO$, the linear form $\phi(u_1,u_2,\bullet)$ is a linear form on $Im\OO$
(which descends to a linear form on $Q=Im\OO/U$), and $\s$ vanishes at $U$ if an
only if this linear form is zero. But for $z\in Im\OO$,
$$\begin{array}{rcl}
\phi(u_1,u_2,z) & = & Re[(u_1z)u_2-(u_2z)u_1] \\
                & = & q((u_1z)u_2,1)-q((u_2z)u_1,1) \\
                & = & -q(u_1z,u_2)+q(u_2z,u_1) \\
                & = & q(z,u_1u_2)-q(u_2,u_1z) \\
                & = & 2q(z,u_1u_2).
\end{array}$$
This is zero for all $u$ if and only if  $u_1u_2=r1$ for some scalar $r$. But multiplying by
$u_1$ on the left, we get $-q(u_1)u_2=ru_1$, thus $r=0$. We conclude that the
zero locus of $\s$ is exactly the set of null-planes in $Im\OO$. (In particular,
it is not empty! Note also that a null-plane must be $q$-isotropic.) This proves our first claim.

\smallskip
Let $d$ be a rank two derivation of $\OO$. Since $d$ has rank two and is skew-symmetric,
we can find two independant vectors $u_1$ and $u_2$ such that
$d(z)=q(u_1,z)u_2-q(u_2,z)u_1.$
Since $d(1)=0$, the plane $U$ generated by $u_1$ and $u_2$ is contained in $Im\OO$.
Since $d$ is a derivation, its kernel $V=U^{\perp}$ is a subalgebra of $\OO$, containing
the unit element. For $v,v'\in V$ and $u\in U$, we get $0=q(u,vv')=q(\bar{v}u,v')$,
hence $V.U\subset U$. This implies that $U$ must be $q$-isotropic, since the
right  multiplication
by a non-isotropic element is invertible. For $u\in U$ nonzero,
consider the right multiplication operator $R_u :\BO\ra \BO$. 
Then $R_u(\BO)$ is a
four dimensional $q$-isotropic subspace of $\OO$.
Since $V$ has codimension $2$ in $\OO$, $V.U$ has codimension at most two in $R_u(\BO)$,
and since it is contained in $U$ we must have $R_u(V)=U$. If $u'\in U$, this means that
we can find $v\in V$ such that $u'=vu$. But then $u'u=(vu)u=-q(u)v=0$. We conclude that
$U$ is a null-plane. Thus
 the projectivization of the space of rank two derivations of $\OO$,
$G(2,Im\OO)\cap\PP\fg_2$, which is   non-empty because $\tdim G(2,Im\OO)=10$   and $\PP\fg_2$ has codimension $7$ in $\PP\fso_7$, can be identified with a subvariety
of $\xadg2$. Being $G_2$-invariant, it must be equal to the adjoint variety.
This proves our second claim.

\smallskip
Our third claim follows. On the one hand, the orthogonal space to a null-plane $U$,
being equal to the kernel of a rank-two derivation, is a six-dimensional subalgebra
of $\OO$. Conversely, we have just proved that the orthogonal to such a subalgebra
is a null plane. \qed

\medskip\noindent
What is the structure of a six dimensional subalgebra $S=U^{\perp}$ of $\OO$?
To understand it, consider another null plane $U_-$, transverse to $S$, and let
$S_-:= U_-^{\perp}$.

\begin{lemma}\label{quat}
$H=S\cap S_-$ is a quaternionic subalgebra.
\end{lemma}

\proof Being the intersection of two subalgebras, $H$ is a subalgebra, and contains $1$.
The hypothesis that $U_-$ be transverse to $S$ is equivalent to the fact that
$H$ is tranverse to $U$ in $S$. In particular, the norm restricts to a nondegenerate
quadratic form on $H$, which must therefore be a quaternionic subalgebra, i.e.,
isomorphic to $\HH$.\qed

\begin{lemma}
The right action of $H$ on $U$ identifies
 $H$ with $\fgl (U)$. \hfill\hfill \linebreak
\indent The scalar product on $\OO$ identifies $U_-$ with $U^*$.
\end{lemma}

We can be even more precise and
explicitly describe   the octonionic multiplication
in terms of the decomposition
$$\OO=\fgl (U)\op U\op U^*.$$
An explicit computation shows that his multiplication is given by the formula
\begin{eqnarray}
(X,u,u^*)(Y,v,v^*)=(XY-2u\ot v^*-2(v\ot u^*)^0,X^0v+Yu,X^tv^*+Y^*u^*),
\end{eqnarray}
and that the norm is
$$q(X,u,u^*)=\det (X)+2\langle u,u^*\rangle.$$
%\smallskip

Here $X^0=\tr(X)I-X$, so that the map $X\mapsto X^0$ is the reflection
in the hyperplane perpendicular to
  the identity. Moreover, $ (X^0)^t$ is the cofactor
matrix of $X$, as $XX^0 =(det X)I$. Also, note that $(XY)^0=Y^0X^0$.
\smallskip

Restricting to $S=\fgl(U)\op U$, we get the multiplication law 
\begin{eqnarray}
(X,u)(Y,v)=(XY,X^0v+Yu),
\end{eqnarray}
while the norm $q(X,u)=\det (X)$ becomes degenerate, with kernel $U$. 

\begin{lemma}\label{decomplemma}
The decomposition $H^{\perp}=U\op U_-$ into the direct sum of two null-planes, is unique.
\end{lemma}

\proof Use the multiplication law $(4)$ on $H^{\perp}$. \qed

\medskip Lemma
\ref{decomplemma} provides   an interesting way to parametrize the set of quaternionic
subalgebras of $\OO$. First note that the Schubert condition $U^{\perp}\cap U_-=0$
defines a $G_2$-invariant divisor $D$ in the linear system $|\cO(1,1)|$ on
$\xadg2\times \xadg2\subset\PP\fg_2\times\PP\fg_2$. This divisor descends to
a very ample divisor in $Sym^2\xadg2$, which we still denote by $D$. Note that
$D$ contains the diagonal, which is the singular locus of $Sym^2\xadg2$.

\begin{proposition}
Let $\cH$ denote the set of quaternionic subalgebras of $\OO$. \hfill\hfill
\linebreak \indent We have the
$G_2$-equivariant identifications
$$\cH\simeq G_2/GL_2\rtimes\ZZ_2\simeq Sym^2\xadg2-D.$$
\end{proposition}

\proof The fact that $\cH$ is $G_2$-homogeneous is well-known, see e.g.
the proof of Theorem 1.27, page 57 of \cite{rosen}.
We prove that the stabilizer $K$ of a quaternionic subalgebra $H$ is isomorphic
to $GL_2$. This stabilizer also preserves $H^{\perp}$, hence, by the
lemma,  the pair $U,U_-$. The subgroup $K^0$ of $K$ preserving $U$, is therefore
either equal to $K$, or a normal subgroup of index two.

For $m\in GL(U)$, the endomorphism $\rho_m$ of $\OO$ defined by
$$\rho_m(X,u,u^*)=(Ad(m)X,mu,{ }^tm^{-1}u^*)$$
is easily checked to be an algebra automorphism, and the map $m\mapsto \rho_m$
defines an isomorphism of $GL(U)$ with $K^0$.

Now let $\s\in K-K^0$. Since the restriction of $\s$ to $H$ is an algebra automorphism,
$$\s(X,u,u^*)=(Ad(s)X,P(u^*),Q(u))$$ for some $s\in GL(U)$, and some invertible
operators $P :U^*\ra U$ and $Q: U\ra U^*$. Composing with a element of $K^0$, we may
suppose that $Ad(s)=1$ and $P\circ Q=\e I$, with $\e=\pm 1$. For $u\in U$ and $u^*\in U^*$,
the condition that $\s(u.u^*)=\s(u).\s(u^*)$ gives $\e =1$ and $\langle Qu,u\rangle =0$.
Thus $Q$ is a skew-symmetric endomorphism from $U$ to $U^*$, and we check that this
is sufficient to ensure that $\s$ is an automorphism. We conclude that $K^0\ne K=K^0
\rtimes\ZZ_2$.

\smallskip
Finally, the decomposition  $H^{\perp}=U\op U_-$, being unique, defines an
injective $G_2$-equivariant morphism from $\cH$ to $Sym^2\xadg2$.
Since $U^{\perp}=H\op U$ and $U_-^{\perp}=H\op U_-$, we have that $U^{\perp}\cap U_-=
U_-^{\perp}\cap U=0$, so that the image of $\cH$ is contained in the complement
of the divisor $D$. In fact there is equality, by Lemma \ref{quat}.\qed

\begin{proposition}
The automorphism group of $S$ fits into an exact sequence
$$1\lra R\lra Aut(S) \lra GL(U)\lra 1,$$
where the radical $R$ is a four-dimensional vector space considered
with its natural abelian group structure. The induced action of $GL(U)$
on $R$ identifies $R$ with $S^3U\ot (\det U)^{-1}$.
\end{proposition}

\proof Recall that the multiplication on $S=H\op U\simeq \fgl (U)\op U$ is given by
$$(X,u)(Y,v)=(XY,X^0v+Yu).$$
An automorphism $\rho$ of $S$ will preserve the norm, hence the kernel $U$ of this
quadratic form. We can therefore write
$$\rho (X,u)=(\rho_2(X),\s(X)+\rho_1 (u)),$$
where $\rho_1\in\GL(U)$, $\s\in Hom(\fgl (U),U)$ and $\rho_2\in Aut(H)$. In particular,
we can find $r\in GL(U)$ such that $\rho_2=Ad(r)$. A straightforward computation shows
that $\rho$ is an automorphism of $S$ if and only if the following conditions hold:
\begin{eqnarray}
\s(XY) =&  \rho_2(X)^0\s(Y)+\rho_2(Y)\s(X) & \forall X,Y\in \fgl (U),\\
\rho_1(Yu) =&  \rho_2(Y)\rho_1(u) & \forall Y\in \fgl (U), u\in U, \\
\rho_1(X^0v) =&  \rho_2(X)^0\rho_1(v) & \forall X\in \fgl (U), v\in U.
\end{eqnarray}
The second condition yields $\rho_1 Y= rYr^{-1}\rho_1$ for all $Y\in \fgl (U)$,
so $r^{-1}\rho_1$ is a homothety and $\rho_2=Ad(\rho_1)$. Since $r^0=(\det r)r^{-1}$,
the third condition follows. Since the first equation is certainly verified by $\s=0$,
the map $\rho\mapsto\rho_1$ defines a surjective morphism from $Aut(S)$ to $GL(U)$.
Note that this surjection is split, since $GL(U)$ can be identified with the subgroup
$Aut_H(S)$ of automorphisms of $S$ preserving $H$.
\smallskip
Consider the kernel of this extension, i.e., the normal subgroup of $Aut(S)$
consisting in morphisms of type $$\rho (X,u)=(X,\s(X)+u),$$
 where $\s\in Hom(\fgl (U),U)$ is subject to the condition that
$$\s(XY) =  X^0\s(Y)+Y\s(X) \qquad \forall X,Y\in \fgl (U).$$
Letting $Y=I$, we see that  $\s(I)=0$. For $Y=X$ we get $\s(X^2)=\tr(X)\s(X)$,
but since $X^2=\tr(X)X-\det(X)I$, this follows from $\s(I)=0$. So the symmetric
part of the condition is fulfilled, and we are left with the skew-symmetric part,
$$\s([X,Y]) =  -2(X\s(Y)-Y\s(X)) \qquad \forall X,Y\in \fsl (U).$$
An explicit computation shows that this defines a four dimensional subspace $R$
of $Hom(\fgl (U),U)$ (take the standard basis $X,Y,H$ of $\fsl(U)$ and check that 
$\s$ is uniquely defined by the choice of $\s(X)$ and $\s(Y)$, which is arbitrary). 
The conjugation action of $\rho\in GL(U)=Aut_H(S)$ is by
$\rho(\s)=\rho^{-1}\circ\s\circ Ad(\rho)$. We finally choose a two dimensional
torus in $GL(U)$ and compute the weights of this action, and they are those of
$S^3U\ot (\det U)^{-1}$. This concludes the proof. \qed

\begin{proposition}\label{prop37}
Let $Aut_S(\OO)$ denote the group of automorphisms of $\OO$  preserving $S$.
The restriction map $Aut_S(\OO)\ra Aut(S)$ is surjective with one dimensional
kernel.
\end{proposition}

Note that $Aut_S(\OO)$ is a maximal parabolic subgroup of $Aut(\OO)=G_2$,
since the adjoint variety $\xadg2=Aut(\OO)/Aut_S(\OO)$.\medskip

\proof We begin with a technical lemma. For $\s\in R\subset Hom(\fgl (U),U)$, let
$\s^{\dagger}\in Hom(U^*,\fgl (U))$ denote its transpose with respect to the
norm on $\OO$. Since the polarisation of the determinant is the symmetric
bilinear form $\det(X,Y)=\tr(X)\tr(Y)-\tr(XY)$, this means that
$$\frac{1}{2}(\tr(Y)\tr\s^{\dagger}(u^*)-\tr(Y\s^{\dagger}(u^*)))=
\langle \s(Y),u^*\rangle \qquad \forall u^*\in U^*, Y\in\fgl (U).$$
Since $\s(I)=0$, we can characterize $\s^{\dagger}$ as the unique morphism from
$U^*$ to $\fsl(U)$ such that
$$\tr(Y\s^{\dagger}(u^*))=-2\langle \s(Y),u^*\rangle
\qquad \forall u^*\in U^*, Y\in\fsl(U).$$

\begin{lemma} For all $\s\in R$, we have the identities
\begin{eqnarray}
\langle\s(u\ot v^*),w^*\rangle = &\langle\s(u\ot w^*),v^*\rangle
&\forall u\in U, v^*,w^*\in U^*\\
 \s^{\dagger}(v^*)u = &-2\s(u\ot v^*) &\forall u\in U, v^*\in U^* \\
 \s^{\dagger}(X^*v^*) = &\s^{\dagger}(v^*)X+2\s(X)\ot v^*&\forall
v^*\in U^*,X\in\fgl(u).
\end{eqnarray}
\end{lemma}

\medskip\noindent {\it Proof of the lemma}. Let $X=u\ot v^*$ and $Y=u\ot w^*$. Then
$XY=\langle u,v^*\rangle Y$ and $X^0=\langle u,v^*\rangle I-u\ot v^*$. The
identity $\s(XY)=X^0\s(Y)+Y\s(X)$ gives
$$\langle u,v^*\rangle \s (u\ot v^*)=(\langle u,v^*\rangle I-u\ot v^*)\s (u\ot w^*)
+(u\ot w^*)\s (u\ot v^*),$$
and the first identity follows. We deduce that for all $w^*\in U^*$,
$$\begin{array}{rcl}
\langle\s^{\dagger}(v^*)u,w^*\rangle & = & \tr(\s^{\dagger}(v^*),u\ot w^*) \\
 & = & -2\langle\s(u\ot w^*),v^*\rangle \\
 & = & -2\langle\s(u\ot v^*),w^*\rangle.
\end{array}$$
This gives the second identity. Finally, for all $Y\in\fgl (U)$, we have
$$\begin{array}{rcl}
\tr (Y\s^{\dagger}(X^*v^*)) & = & -2\langle \s(Y), X^*v^*\rangle \\
   & = & -2\langle X^0\s(Y), v^*\rangle \\
   & = & -2\langle \s(XY), v^*\rangle +2\langle Y\s(X), v^*\rangle \\
   & = &  \tr (XYs^{\dagger}(v^*))+2\tr Y(s(X)\ot v^*),
\end{array}$$
and this implies the last identity. \qed

\begin{lemma}
For all $\s\in R$, the map $d_{\s}\in End(\OO)$ defined by
$$d_{\s}(X,u,u^*)=(-\s^{\dagger}(u^*),\s(X),0), \qquad X\in\fgl (U),
u\in U, u^*\in U^*,$$
is a derivation of $\OO$.
\end{lemma}

\proof Easy verification with the formulas of the previous lemma. \qed

\begin{lemma}
For all $\rho\in \fgl (U)$, the map $d_{\rho}\in End(\OO)$ defined by
$$d_{\rho}(X,u,u^*)=(ad(\rho)X,\rho(u),-\rho^t(u^*)),
\qquad X\in\fgl (U), u\in U, u^*\in U^*,$$
is a derivation of $\OO$.
\end{lemma}

\proof Straightforward. \qed

\medskip
We can now complete the proof of   Proposition \ref{prop37}. The differential at the
identity of the map $Aut_S(\OO)\ra Aut(S)$ is the natural restriction map
$Der_S(\OO)\ra Der(S)$. The two previous lemmas imply that this map is surjective.
The map $Aut_S(\OO)\ra Aut(S)$ is therefore surjective as well.

Consider some automorphism $\gamma$ of $\OO$ acting trivially on $S$.
The simple fact that it preserves the norm implies that
$$\gamma (X,u,u^*)=(X,u+\delta (u^*),u^*)$$
for some skew-symmetric map $\delta : U^*\ra U$. Up to scale, there is
only one such skew-symmetric map. Moreover, being a rank two skew-symmetric
map with a null plane for image, it must be a derivation. We conclude that
the kernel of the restriction map $Aut_S(\OO)\ra Aut(S)$ is the additive
group of automorphisms of the form $I+d$, $d$ a rank two derivation with image
in $U$. \qed

\medskip

\begin{definition}
In what follows, we fix a six dimensional subalgebra $S$ of $\BO$,
denote it by $\BS$ 
and call it the {\it sextonion algebra}.
Recall formula $(5)$, which gives a model of the sextonion algebra over an 
arbitrary field, for example over the real numbers :
 $\SS\simeq\fgl(U)\op U$ for some two-dimensional vector space $U$, 
and the product is given by the simple formula
$$(X,u)(Y,v)=(XY,X^0v+Yu),$$
where $X^0=\tr(X)I-X$. We get a six-dimensional alternative algebra,
with zero divisors. 
\end{definition}

\section{Review of the triality and $r$-ality constructions}\label{ralitysect}

For $\AA$ a composition algebra, define the {\it triality group}
$$T(\AA)=\{\th=(\th_1,\th_2,\th_3)\in SO(\AA)^3\mid \;\;
\th_3(xy)=\th_1(x)\th_2(y)\;\forall x,y\in\AA\}.$$
There are three natural actions of
$T(\AA)$ on $\AA$ corresponding to its three projections on $SO(\AA)$,
and we denote these representations by $\AA_1$, $\AA_2$, $\AA_3$.
See \cite{LMadv} for more details.
We let $\ft(\BA)$ denote the corresponding Lie algebra.

Now let $\AA$ and $\BB$ be two composition algebras. Then
$$\fg(\AA,\BB) = \ft(\AA)\times\ft(\BB)\oplus (\AA_1\ot\BB_1)
\oplus (\AA_2\ot\BB_2) \oplus (\AA_3\ot\BB_3)
$$
is naturally a semi-simple Lie algebra when $\BA,\BB$ are among
$0,\BR,\BC,\BH,\BO$.

The triality Lie algebras can be generalized to $r$-ality for all $r$
to recover the generalized Freudenthal chart (see \cite{LMmagic}). For $r>3$ we have
$$
\ft_r(\RR)=0, \qquad\ft_r(\CC)=\CC^{\op (r-1)},
\qquad\ft_r(\HH)=\fsl_2^{\times r}\qquad\ft_r(\BS)=\fsl_2^{\times r}\op
\BC^{2(r-1)}
$$
and
$$\fg_r(\AA,\BB) =  \ft_r(\AA)\times\ft_r(\BB) \op
\bigoplus_{1\leq i<j\leq r}\AA_{ij}\ot\BB_{ij}.$$
The algebras $\fg_r(\AA,\BB)$ are all semi-simple when
 $\BA,\BB$ are among
$0,\BR,\BC,\BH,\BO$, and moreover, all simple Lie algebras
except $\fg_2$ arise by this construction ($\fg_2$ can
be recovered by supplementing this list with the derivation algebras).

\medskip
The goal of this section is to show that this construction works with the sextonions,
which is not a complexified composition algebra since its natural quadratic form
is degenerate. Nevertheless, the definitions of the triality group and algebra make
sense.

\begin{proposition}
The triality algebra $\ft(\SS)=Der(\SS)\op Im(\SS)^{\op 2}$. Its dimension is $18$.
\end{proposition}

\proof Same proof as in Barton \& Sudbery,  \cite{bs}. \qed

\medskip
There is no natural inclusion of $\ft(\SS)$ in $\ft(\OO)$, but the subalgebra
$\ft_{\SS}(\OO)\subset \ft(\OO)$ of triples $\th\in\fso(\OO)$ such that
$\th_i(\SS)\subset\SS$ for $i=1,2,3$, is a kind of substitute for $\ft(\SS)$.

\begin{corollary}
The natural morphism $\ft_{\SS}(\OO)\ra \ft(\SS)$ is surjective
with one dimensional kernel.
\end{corollary}

This allows one to determine the structure of $\ft(\SS)$. Indeed,
$\ft_{\SS}(\OO)$ is the subalgebra of $\ft(\OO)$ preserving $\SS$. This is the
same as preserving its orthogonal, which is a null plane $U$. The Grassmannian
of isotropic planes in $\OO$ is homogeneous under the action of $T(\OO)=Spin_8$;
in fact, it is
  the adjoint variety of $Spin_8$. The stabilizer
of an isotropic two plane, for example the stabilizer of $U$, is therefore a
maximal parabolic subgroup, which can also be defined as the stabilizer of a
highest root space.
Recall that the choice of a highest  root space in $\ug=\fso_8$ induces a 5-grading
$\ug=\ug_{-2}\op\ug_{-1}\op\ug_0\op\ug_1\op\ug_2$. The stabilizer of the highest
root space $\ug_2$ is $\ug_0\op\ug_1\op\ug_2=\ft_{\SS}(\OO)$. We let
$\ft^*(\SS)=[\ug_0,\ug_0]\op\ug_1\op\ug_2$, the   intermediate subalgebra
of $\ug =\fso_8$.

\smallskip
Using \cite{LMpop}, section 3.4, we can prove:
\begin{proposition}
If $A,B,C,D$ are 2-dimensional vector spaces, we have identifications

$$\begin{array}{ccccc}
 \ft(\OO)& = & \fsl(A)\times\fsl(B)\times\fsl(C)\times\fsl(D) &\op & A\ot B\ot C\ot D \\
    \cup & &       \cup & & \cup \\
 \ft^*(\SS)& = & \fsl(A)\times\fsl(B)\times\fsl(C)\op\CC &\op & A\ot B\ot C  \\
    \cup & &       \cup & & \\
 \ft(\HH)& = & \fsl(A)\times\fsl(B)\times\fsl(C) & &
\end{array}$$
\end{proposition}

The triality algebra $ \ft(\SS)$ is then also identified with
$\fsl(A)\times\fsl(B)\times\fsl(C)\op\CC
\op  A\ot B\ot C$, but the Lie algebra
structure is not exactly the same as that of  $ \ft^*(\SS)$.
(In the notation of the introduction, $\ft (\BS)$ corresponds to the
intermediate algebra $\fg ''$.)
\bigskip

We include the sextonions in the triality construction by letting
\begin{eqnarray}
\nonumber
\fg(\AA,\SS)^+ &= &\ft(\AA)\times\ft_{\SS}(\OO)\oplus (\AA_1\ot\SS_1)
\oplus (\AA_2\ot\SS_2) \oplus (\AA_3\ot\SS_3), \\
\nonumber
\fg(\AA,\SS) &= &\ft(\AA)\times\ft^*(\SS)\oplus (\AA_1\ot\SS_1)
\oplus (\AA_2\ot\SS_2) \oplus (\AA_3\ot\SS_3).
\end{eqnarray}
We can also define $\fg(\SS,\SS)$ by replacing $\AA$ with $\SS$ and $\ft(\AA)$ with
$\ft^*(\SS)$ in this formula. Since $\ft^*(\SS)$ is a subalgebra of $ \ft(\OO)$,
$\fg(\AA,\SS)$ is defined as a subvector space of $\fg(\AA,\OO)$.

\begin{proposition}
$\fg(\AA,\SS)$ is a Lie subalgebra of $\fg(\AA,\OO)$.
\end{proposition}

\proof From the definition of the Lie bracket of $\fg(\AA,\OO)$ given in
\cite{LMadv}, we see that we just need to check that the maps $\Psi_i :
\La 2\OO_i\ra\ft(\OO)$ take $\La 2\SS_i\subset\La 2\OO_i$ inside
$\ft^*(\SS)\subset\ft(\OO)$. This is clear for $\Psi_1$, since the
image of $\Psi_1(s,s')$ is just the plane generated by $s$ and $s'$.
This is also clear for $\Psi_2$ and $\Psi_3$:
$\Psi_2(s,s')$ and  $\Psi_3(s,s')$ are defined in terms of left and right multiplication
by $s$ or $s'$, so that the subalgebra $\SS$ is preserved when $s$ and $s'$ belong to it.
\qed

\begin{proposition}
For $\AA\ne\SS$, $\fg(\AA,\SS)$ is the intermediate subalgebra of the simple
Lie algebra $\fg(\AA,\OO)$.    
It is a maximal parabolic subalgebra minus the one dimensional center,
and its semi-simple part is equal to the simple Lie algebra $\fg(\AA,\HH)$.
\end{proposition}

\proof We saw in \cite{LMadv} that a Cartan subalgebra of $\fg(\AA,\OO)$ is given
by the product of two Cartan subalgebras in $\ft(\AA)$ and $\ft(\OO)$. Moreover,
once we have chosen a set of positive roots for $\ft(\OO)=\fso_8$, its highest root
can be chosen as a highest root for $\fg(\AA,\OO)$. Since $\fso_8=\La 2\OO$, we can
identify the highest root line with an isotropic two plane in $\OO$, which we can
choose to be the null-plane $\BH^{\perp}$. It is then straightforward to check
that the stabilizer of the highest root line in the adjoint representation $\fg(\AA,\OO)$
is exactly $\fg(\AA,\SS)^+$, and our first claim follows. The second claim is a simple exercise.
Note that  $\fg(\AA,\HH)$ is embedded in $\fg(\AA,\OO)$ through the natural
embedding of $\ft(\HH)\simeq\ft_{\HH}(\OO)\subset\ft(\OO)$, as explained in
\cite{LMpop}, section 3.6. \qed

\begin{definition}
We denote by $\fe_{7\frac 12}$ the algebra $\fg(\SS,\OO)$, which is 
intermediate between the exceptional algebras 
$\fe_7=\fg(\HH,\OO)$ and $\fe_8=\fg(\OO,\OO)$.
\end{definition}

%\bigskip
Although we have no direct proof, we observe that, for $a,b=0,1,2,4,6,8$:
\begin{align}
\tdim Der(\BA) &= \frac{4 (a-1) (a-2)}{a+4}
\\
\tdim \ft (\BA) & = \frac{6a (a-1)  }{a+4}
\\
\tdim\fg(\BA,\BB) &=
 \frac{ 3 (4 a + a b + 4 b - 4) (2 a + a b + 2 b)}
 {(a + 4) (4 + b) }
\end{align}

\bigskip
Here are the resulting algebras $\fg(\BA,\BB)$ giving
rise to an expanded magic chart. The first row is the dimension of $\BA$. The first column
contains the derivation algebras:

\begin{center}
\begin{tabular}{|c|c|c|c|c|c|c|} \hline
$-2/3$ & $0$ & $1$ & $2$ & $4$ & $6$ & $8$ \\ \hline\hline
$0$ & $0$ & $A_1$ & $A_2$ & $C_3$ & $C_3.H_{14}$ & $F_4$ \\ \hline
$0$ & $T_2$ & $A_2$ & $2A_2$ & $A_5$ & $A_5.H_{20}$ & $E_6$ \\ \hline
$A_1$ & $3A_1$ & $C_3$ & $A_5$ & $D_6$ & $D_6.H_{32}$ & $E_7$ \\ \hline
$A_1.H_4$ & $(3A_1).H_8$ & $C_3.H_{14}$ & $A_5.H_{20}$%
 & $D_6.H_{32}$ & $D_6.H_{32}.H_{44}$ & $E_7.H_{56}$ \\ \hline
$G_2$ & $D_4$ & $F_4$ & $E_6$ & $E_7$ & $E_7.H_{56}$ & $E_8$ \\ \hline
\end{tabular}
\end{center}

\medskip
The convention here is that a Lie algebra $G.H_{2n}$ means that the
Lie algebra of type $G$ has a representation $V$ of dimension $2n$
which admits an invariant symplectic form $\omega$. Then $G$ acts on the
Heisenberg algebra of $(V,\omega)$ and $G.H_{2n}$ denotes the semi-direct
product. These algebras are not reductive and the Heisenberg algebra
is the radical.

\medskip
There is another series of Lie algebras, the
Barton-Sudbery intermediate Lie algebras
of \cite{bs}. These are called intermediate because they are
intermediate between the derivation algebras and the triality algebras.
This gives the following table:

\begin{center}
\begin{tabular}{|c|c|c|} \hline
$0$ & $0$ & $0$ \\ \hline
$0$ & $T_1$ & $T_2$ \\ \hline
$A_1$ & $2A_1$ & $3A_1$ \\ \hline
$A_1.H_4$ & $2A_1.H_6$ & $3A_1.H_8$ \\ \hline
$G_2$ & $B_3$ & $D_4$ \\ \hline
\end{tabular}
\end{center}

%\medskip

\section{Universal decompositions}\label{univdecompsect}

Let $\fg $
be an intermediate Lie algebra and write $V=\ug_1$.

\subsection{Decomposition of $\fg\ot\fg$}
In order to decompose $S^2\fg$, $\La 2\fg$, we need to
understand the decomposition of $\fh\ot V$. This turns out
to be   uniform:

\begin{proposition} Let $\fg= \fh\op V\op \BC$ be an intermediate
Lie algebra. Then
$$
\fh\ot V= \fh V\op V\op (\fh V)_{Aad}
$$
where $(\fh V)_{Aad}$ is  as follows:
$$
\begin{matrix}
\fh & V & (\fh V)_{Aad} \\
\fsl_2\times\fso_n &  \BC^2\ot\BC^n = W\ot V_{\o_1} & W\ot (V_{\o_1}\op V_{\o_3})\\
\fsl_n & \BC^n\op\BC^{n*}=V_{\o_1}\op V_{\o_{n-1}} &
V_{\o_1+\o_{n-2}}\op V_{\o_2+\o_{n-1}}\\
\fc_n & V_{\o_1} & V_{\o_1+\o_2}
\end{matrix}
$$
and from \cite{LMseries} we recall for the subexceptional series:
$$\begin{array}{lcccccc} %\hline
 & A_1 & A_1^{\op 3} & C_3 & A_5 & D_6 & E_7 \\ %\hline
  & & & & & & \\
V & \left[3\right] & \left[1,1,1\right] & \left[0,0,1\right] &
\left[0,0,1,0,0\right]
& \left[0,0,0,0,0,1\right] & \left[0,0,0,0,0,0,1\right] \\
\fh & \left[2\right] & \left[2,0,0\right] & \left[2,0,0\right] &
\left[1,0,0,0,1\right]
& \left[0,1,0,0,0,0\right] & \left[1,0,0,0,0,0,0\right] \\
(\fh V)_{Aad} & \left[1\right] & \left[1,1,1\right]\ot\rho & \left[1,1,0\right] &
\left[1,1,0,0,0\right]
& \left[1,0,0,0,1,0\right] & \left[0,1,0,0,0,0,0\right]\\
 \end{array}$$

\bigskip In the column corresponding to $A_1^{\op 3}$,
$\rho$ denotes the two-dimensional irreducible representation of
$\Gamma={\mathfrak S}_3$.
\end{proposition}

Given two modules $V,W$, the module
$(VW)_{Aad}$ is defined and discussed in \cite{LMseries}, section 2.3.
In particular, in most cases it may be determined by pictorial methods using
Dynkin diagrams. In the case $W\subseteq I_2(X)\subset S^2V^*$, where $X$ is the
closed $G$ orbit in $\BP V$, then $(VW)_{Aad}$ is a space of linear
syzygies among the quadrics in $W$.

Recall the universal decomposition formulas
of Vogel $\La 2\fh= \fh \op \fh_{2}$, $S^2\fh
=\fh^2\op \fh_{Q}\op\fh_{Q'}\op \BC$ \cite{vog2}.
We obtain uniform
decompositions of
\begin{eqnarray}
\nonumber S^2\fg & = &S^2\fh\op S^2V\op \fh\ot V
\op\fh\op V\op \BC, \\
\nonumber \La 2\fg &= &\La 2\fh\op \La 2V\op \fh\ot V
\op\fh\op V.
\end{eqnarray}
Here, Vogel's decompositions work if we take
\begin{eqnarray}
\nonumber\fg_2 &= &\fh_2\op V_2\op \fh V\op (\fh V)_{Aad} \op \fh \op V, \\
\nonumber \fg^2 &= &\fh^2\op\fh V\op V^2\op\fh\op V\op\CC, \\
\nonumber \fg_Q &= &\fh_Q\op (\fh V)_{Aad} \op V\op \fh
\end{eqnarray}
(the last equation assumes we are in the case $\fh_{Q'}=0$).
 It would be interesting to determine to what extent the
Cartan powers of $\fg_Q$ satisfy the dimension formulas of \cite{LMuniv}.

\subsection{Cartan powers of $\fg$}
One can check that the formulas above really define $\fg$-submodules of $\fg\ot\fg$.
For example, $\fg_2$ is the $\fg$-submodule of $\La 2\fg$ generated by $\fh_2$. 

In general, given an irreducible finite dimensional
 $\ug$-module $W_{\l}$ with highest weight line $\ell_{\l}$, 
we can define a  highest weight $\fg$-module
module $V_{\l}$ by taking   $V_{\l}=U(\fg)\ell_{\l}$. Note that
this is the same as taking $V_{\l}=U(\ug_1)W'_{\l}=S^{\bullet}(\ug_1)W'_{\l}$, where
$W'_{\l}$ is the $\fh$-module $U(\fh ).\ell_{\l}$. As in \cite{shtepin}, where the 
case of classical intermediate algebras was studied, weights
$\l,\mu$ of $\ug$ will give rise to the same $\fg$ module if and only if they project
to the same weight in the weight lattice of $\fh$ (considered as a subspace
of the weight lattice of $\ug$).

In general we have no effective way of computing $V_{\l}$ from $W_{\l}$
but we do have the following special case:

\begin{proposition}\label{prop51} Suppose that the highest weights of $V$ and $\fh$ are 
linearly independant. Then, as an $\fh$-module,
$$\fg^{(k)}=\bigoplus_{p+q\le k}\fh^{(p)}V^{(q)}.$$
\end{proposition}

\proof As a subspace of $S^k\fg$,  the Cartan power $\fg^{(k)}$
is generated by the powers $x^k$ of the highest weight vectors of
$\fh$, and their images by successive applications of vectors in
$V=\fg_1$. For $v,w\in V$, we have
$$\begin{array}{rcl}
ad(v)x^k & = & kx^{k-1}(xv), \\
ad(w)ad(v)x^k & = &  k(k-1)x^{k-2}(xv)(xw)+kx^{k-1}\omega(xv,w).
\end{array}$$
First observe that the last expression is symmetric in $v$ and $w$,
so that the action of $V$ induces an action of $Sym(V)$.
Second, the last term is a multiple of $x^{k-1}$, and that kind
of terms generate $\fg^{(k-1)}\fg_2$. By induction on $k$,
we are reduced to proving that the $\fh$-module spanned by
tensors of the form $x^{k-q}(xv)^q$, for $x$ a highest weight
vector of $\fg_0$ and $v\in V$, is a copy of
$\fh^{(k-q)}V^{(q)}$. Since it follows from the hypothesis that 
the weights of these modules, as 
$k$ and $q$ vary, are distinct, Schur's lemma will imply our claim. 

We first prove that we can suppose
that $xv$ is a highest weight vector in $V$. To see this, recall
that the image of $x$ in $\PP V$
is the linear space denoted $\langle q^x\rangle$ in \S 2, and
$\langle q^x\rangle\cap X^{ad}_H$ is a 
smooth quadric hypersurface in $\langle q^x\rangle$ whose equation is
$$0=q^x(xu,xv)=\omega(u,xv)=\omega(v,xu).$$
The first expression shows that this quantity does not depend
on $v$, but only on $xv$, and the second one shows that it does
not depend on $u$, but only on $xu$. Let $V_x=x.V\subset V$ denote
the deprojectivization of $\langle q^x\rangle$. We get
$S^qV_x=S^{(q)}V_x\op q^xS^{(q-2)}V_x\op\cdots$. Since $q^x$ is given
by expressions of type $\omega(v,xu)$, it must be considered as belonging
to $\fg_2$, and we remain with $S^{(q)}V_x$ only. By definition,
this space is generated by $q$-th powers of vectors that belong to the
quadric hypersurface $q^x=0$, hence also to the cone over the
closed $G_0$-orbit in $\PP V$. This proves our claim that we can
suppose $xv$ to be a highest weight vector.

The stabilizers $\fg_0^x$ and $\fg_0^{xv}$ are two parabolic
subalgebras of $\fg_0$. Their intersection must therefore contain
a Cartan subalgebra, and we can choose a Borel subalgebra containing
this Cartan subalgebra and contained in $\fg_0^x$. In other words,
we may suppose that $x=x_{\ta}$ is a highest root vector of $\fg_0$,
while $xv$ is a weight vector of $V$. Of course we can also suppose
that $v$ itself is a weight vector, say of weight $\mu$, so that
the weight of $xv$ is $\mu+\ta$.

Now we use the fact that $V$ is a minuscule $\fg$-module. In particular,
$\mu(H_{\ta})$ and $(\mu+\ta)(H_{\ta})=\mu(H_{\ta})+2$ belong to
$\{-1,0,+1\}$, hence $\mu(H_{\ta})=-1$. For simplicity, suppose that
$\fg_0$ is not of type $A$, so that the highest root is a multiple
of a fundamental weight $\o_{\gamma}$. Then $(\mu+\ta)(H_{\g})
>\mu(H_{\g})\ge -1$, so $(\mu+\ta)(H_{\g})\ge 0$. If $(\mu+\ta)(H_{\b})\ge 0$
for every simple root $\b$, then $\mu+\ta$ is a dominant weight, hence
the highest weight of $V$.  If $(\mu+\ta)(H_{\b})<0$ for some
simple root $\b$, necessarily distinct from $\g$, then the corresponding
reflection stabilizes $\ta$, but changes $\mu+\ta$ into the greater
root $\mu+\ta+\b$. By induction, we may therefore suppose that
the weight $\mu+\ta$ of $xv$ is the highest weight of $V$.
Then $x^{k-q}(xv)^q$ is a highest weight vector of the Cartan
product $\fg_0^{(k-q)}\fg_1^{(q)}$, and we are done. \qed

\rem The only case of rank greater than two, for which the hypothesis of   
Proposition \ref{prop51} does not hold, is when $\fg=\fsp_{2n+1}$ is an odd symplectic
Lie algebra. The highest weight of $\fh=\fsp_{2n}$ is $2\o_1$, twice the
weight of $V=\CC^{2n}$. The Proposition does not hold in that case, but 
it is easy to see that as an $\fsp_{2n}$-module
$$\fsp_{2n+1}^{(k)}=S^{2k}\CC^{2n}\op S^{2k-1}\CC^{2n}\op\cdots\op\CC^{2n}\op\CC
\simeq S^{2k}(\CC^{2n}\op\CC).$$

\rem Consider the intermediate Lie algebra of $\fsl_{n+2}$, that we denote by
$${\widetilde{\fsl}}_{n+1}=\fsl_n\op\CC^n\op (\CC^n)^*\op\CC.$$
By the previous theorem, the decomposition of its Cartan powers into 
$\fsl_n$-modules is 
$${\widetilde{\fsl}}_{n+1}^{(k)}=\bigoplus_{p+q+r\le k}V_{(p+q)\o_1+(p+r)\o_{n-1}}.$$
This is exactly the formula for the restriction of the $\fsl_{n+1}$-module $\fsl_{n+1}^{(k)}$ 
to $\fsl_n$ given by the usual branching rule. We therefore have two  
different Lie algebras, $\fsl_{n+1}$ and ${\widetilde{\fsl}}_{n+1}$, not only with the 
same dimension, but such that in any degree, their Cartan powers have the same dimensions.

\section{Decomposition formulas in the magic chart}\label{rowbyrowsect}

The sextonions allow one to add a new column to Freudenthal's magic square. 
We know that for each row of the original square, there are a few prefered 
representations, leading to nice dimension and decomposition formulas 
for some of their plethysms (see \cite{LMseries}). In this section we 
address the problem of extending these results to the sextonionic case. 
What the prefered representations should be is easy to imagine: take 
a prefered representation $V_{\OO}$ from the octonionic column; it contains a 
prefered representation $V_{\HH}$ from the quaternionic column, and the sextonionic
representation $V_{\SS}$ is simply the $\fg$-submodule  of $V_{\OO}$ generated by
$V_{\HH}$, where $\fg$ is the intermediate Lie algebra.

We adopt the notation $V_0= \BC\op \ug_1 $. In several of the
 modules below $V_0$ will replace the trivial representation
 in the decomposition formulas. This makes sense when the trivial
 representation corresponds to the copy of $\BC$ in $\ug_0$.

\subsection{First row}
Here  we have one distinguished representation, call it
$V= \cJ_3(\BS)_0$. It is the complement of the
symplectic form in $\La 2\BC^6\op \BC^{6*}=\La 2(\BC^6+\BC)$.

As   graded $\fsp_6$-modules, we have
\begin{align}
V&=V_{\o_2} \op V_{\o_1}\\
\fg &=V_{2\o_1 } \op V_{\o_3} \op \BC\\
V_2 &= V_{\o_1+\o_3}\op (V_{\o_1+\o_2}\op V_{\o_1}) \op V_{\o_2}
\end{align}

We have the following decomposition formulas, which
agree with those in \cite{LMseries}:
\begin{eqnarray}
S^2V &= &V^2\op V\op V_0\\
\La 2 V &= &\fg\op V_2 .
\end{eqnarray}

\subsection{Second row}
Here, we have two dual distinguished representations, call one of them
$V=\cJ_3(\BS)= \La 2(\BC^6\op\BC)$.
 As   graded $\fsl_6$-modules, we have
\begin{align}
V&=V_{\o_2} \op V_{\o_5}\\
V^*&=V_{\o_4} \op V_{\o_1}\\
\fg &=V_{\o_1+\o_5} \op V_{\o_3} \op \BC\\
%V_2 &= V_{\o_1+\o_3}\op (V_{\o_3+\o_4}\op V_{\o_1}) \op V_{\o_4}\\
V_0&=\BC\op V_{\o_3}
\end{align}
The $\fsl_6$-module $V$ is exceptional in the sense of \cite{brion}
and its symmetric algebra behaves the same as the rest of the Severi series,
namely
$$
S^dV=\bigoplus_{i+2j+3k=d}V^{(i)}(V^*)^{(j)}(V_0)^{(k)}
$$
where we take $\fsl_6$-Cartan products in the factors.

\subsection{Third row}\label{thirdrowsubsect}
There are three distinguished representations, which we
call $V=\cZ_2(\BS), V_2,\fg$. As graded $\fso_{12}$-modules they are
\begin{align}
V&=V_{\o_6} \op V_{\o_1}\\
\fg &=V_{\o_2} \op V_{\o_5} \op \BC\\
V_2 &= V_{\o_4}\op (V_{\o_1+\o_6}\op V_{\o_5})  \op V_{\o_2}\\
V_0&= \BC \op V_{\o_5}
\end{align}
Here again, $V$ is exceptional in the sense of \cite{brion}
and its symmetric algebra behaves the same as the rest of the
subexceptional series, namely
$$
S^dV=\bigoplus_{i+2j+3k+4l+4m=d}
V^{(i+k)}\fg^{(j)}V_0^{(l)}V_2^{(m)}
$$

Here some care must be taken in interpreting the formula.
In Brion's list there are $14$ generators of the symmetric algebra which
do not co\"{\i}ncide with the generators we use.
The critical difference is that the   product $\fg V$, is not the
Cartan product as $\fsl_6$-modules, but instead
$$
\fg V= V_{\o_2 +\o_6} \op (V_{\o_5 +\o_6}\op V_{\o_3})\op V_{\o_5}
$$
where note that the $V_{\o_3}$ would not appear in the $\fsl_6$ Cartan
product. All other   products coincide
with the Cartan product in $\fsl_6$.
Thus the interpretation of the algebra structure is different.

\smallskip
 The justification for $\fg V$ is as follows.
 In $V_{\o_2}\ot V_{\o_6}$,
 the submodule $V_{\o_2+\o_6}$ is generated by tensors of the form $P\ot S$
 with $P\in G_Q(2,12)$,   $S\in \BS_6\subset G_Q(6,12)$ an isotropic
 $6$-plane, where $P\subset S$.  We have a map
 $V_{\o_5}\ot V_{\o_6}\ra V_{\o_1}$, which may be seen geometrically
 as follows. Let $S'\in \BS_6'\subset \BP V_{\o_5}$ be a $6$-plane in the
 other  family. Generically $S\cap S'$ is a point of the quadric. This
 geometric intersection extends to a linear map
 $V_{\o_5}\ot V_{\o_6}\ra V_{\o_1}$. The action of $V_{\o_5}$ thus
 produces tensors of the form
 $v\ot P$ with no incidence condition on $v$ and $P$, in particular
 a projection to $V_{\o_3}$ by wedging them together. 

%In addition  we must have everything from the $\fe_7$ Cartan product and now
% a dimension count shows that we cannot have more.

\section{Dimension formulas}

We have the following generalizations of the theorems
in \cite{LMadv}:

\begin{theorem}\label{dimexc}
Let $\fg =\fsl_2, \fsl_3, \fg_2, \fso_8, \ff_4, \fe_6,
 \fe_7, \fe_{7\frac{1}{2}}=\fe_7\op V_{\o_7}\op\BC, \fe_8$,
with respectively $a=-4/3,-1,-2/3, 0,1,2,4,6, 8$. Then for all $k\ge 0$,
$$\dim \fg^{(k)} = \frac{3a+2k+5}{3a+5} \frac{\binom{k+2a+3}{k}
\binom{k+\frac{5a}{2}+3}{k}\binom{k+3a+4}{k}}
{\binom{k+\frac{a}{2}+1}{k}\binom{k+a+1}{k}}.$$
\end{theorem}

\begin{theorem}
Let $V$ be the distinguished module, of dimension $6a+8$, of
a Lie algebra $\fg$ in the subexceptional series,
with $a=-\frac{2}{3}, 0, 1, 2, 4,6, 8$. Then
$$\begin{array}{rcl}
\dim \fg^{(k)} & = & \frac{2k+2a+1}{2a+1}\frac{
\binom{k+\frac{3a}{2}-1}{k}\binom{k+\frac{3a}{2}+1}{k}
\binom{k+2a}{k}}{\binom{k+\frac{a}{2}-1}{k}
\binom{k+\frac{a}{2}+1}{k}}, \\
\dim V^{(k)} & = & \frac{a+k+1}{a+1}\frac{\binom{k+2a+1}{k}
\binom{k+\frac{3a}{2}+1}{k}}{\binom{k+\frac{a}{2}}{k}}, \\
\dim  V_2^{(k)} & = & \frac{ ( 4k+3a+2)}{( k+1)(3a+2)}
\frac{\binom{k+a}{k}\binom{ k+a+1}{k}\binom{k+\frac{3a}{2}-1}{k}
\binom{ k+\frac{3a}{2}}{k}\binom{ 2k+2a+1}{2k}}
{\binom{k+\frac{a}{2}-1}{k}\binom{ k+\frac{a}{2}}{k}
\binom{ 2k+a}{2k}}.
\end{array}$$
\end{theorem}

\begin{theorem} Let $V$ be the distinguished module
in the Severi series, with
$a=-\frac{2}{3}, 0, 1, 2, 4, 6, 8$. Then
$$
{\rm dim}\;V^{(k)}   =
  \frac{(2k+a)(k+ a)}{a^2}
\frac{\binom{k+a-1}{k}\binom{k +\frac{3a}{2}-1}{k }
 }{\binom{k +\frac{a}{2}}{k}}.
 $$
\end{theorem}

Unfortunately our proofs are just case by case applications of
the Weyl dimension formulas, plus the decomposition formulas from 
Proposition 5.2 and the previous section. Even then we obtain
in each case a polynomial $P(k)$ of
the correct degree, but that is not obviously the same polynomial
as obtained above. To check we used Maple to test that the two polynomials
agree on $\tdeg P+1$ points and therefore must be equal.

Here are outlines of the proofs:

 \subsection*{Severi case}  The triality formula for $a=6$ predicts
$$
\tdim V^{(k)}=
 \frac{(2k+6) (k+6)}{36}\frac{ \binom{k+5}{k}\binom{k+8}{k}}
 {\binom{k+3}{k}}
 $$
which, as a function of $k$ is a polynomial of degree $12$.
We compare with the Weyl dimension formula
applied to the $\fsl_6$-module
$$
V^k=(V_{\o_2}\op V_{\o_5})^{(k)}= \sum_{i=0}^k V_{(k-i)\o_2+ i\o_5}
$$
which also gives a polynomial of degree $12$ in $k$.
First note that for all positive roots $\a$ we
have $(\o_2,\a)$ and $(\o_5,\a)$ either $0$ or $1$. Seperate the positive
roots of $\fsl_6$ into four groups accordingly:
$\Delta_{0,0},\Delta_{1,0},\Delta_{0,1},\Delta_{1,1}$ where
the first subscript is $(\o_2,\a)$ and the second is $(\o_5,\a)$.
$\Delta_{0,0}$ has four elements, three of which have
$(\rho,\a)=1$ and one with $(\rho,\a)=2$.
$\Delta_{1,0}$ has three elements, with $(\rho,\a)=1,2,3$.
$\Delta_{0,1}$ has six elements, with $(\rho,\a)=1,2,2,3,3,4$.
$\Delta_{0,1}$ has two elements, with $(\rho,\a)=4,5$.
Thus the numerator in the WDF becomes
$$
2(k+5)(k+4)  [ \sum_{i=0}^k
(k-i +2)^2 (k-i +3)^2 (k-i +1) (k-i +4) (i +3) (i +2) (i +1)  ]
$$
Dividing by the denominator, and considering, e.g., the
$i=[k/2]$ term), we obtain another polynomial that is a sum of $k+1$
terms of
degree $11$ in $k$, but these terms collapse by using formulas
for $\sum_k i^k$ to give a polynomial of degree twelve.
  One then easily checks they agree for the first
$13$ values of $k$ so they must be equal.

 \medskip

\subsection*{Exceptional row}
For $\fg^k$ in the  exceptional row the dimension of the
relevant $\fe_7$ modules are as follows:
$$
\dim V_{i\o_1+j\o_7}  = \frac{
(j+5)(2i/17 +j/17 +1)\binom {j+9}{ 9}\binom{11+i}{ 11}\binom{8+i}{8}\binom{16+i+j}{16}\binom{13+i+j}{13}
}
{ 5\binom{3+i}{3}\binom{8+i+j}{8}\binom{5+i+j}{5}}
$$
For $\tdim \fg\up k$, one takes the sum over $i+j\le k$
and compares it with the triality formula. Both are polynomials
of degree $45$ in $k$ but they are not obviously equal so we
evaluated them both at $45$ points (plus zero) to check their 
equality.

\medskip

\subsection*{Subexceptional row}
We calculate as above.
The relevant dimensions of the
 $\fso_{12}$-modules
that
need to be summed over are respectively
\begin{align*}
\tdim V_{i\o_5+j\o_2}
&= {\scriptstyle 
\frac{
(2i+j+9)(j+3)\binom{i+5}{5}\binom{i+4}{4}\binom{i+j+8}{8}
\binom{i+j+7} {7}\binom{j+5}{5}
}
{27(i+1)\binom{i+j+3}{3}\binom{i+j+4}{4}
 }}
\\
\tdim V_{a\o_4+b(\o_1+\o_6)+ c\o_5+d\o_2}
&=
{\scriptstyle \frac{1}{158018273280000}(1+b)^2 (2+b+d) (3+b+d) (4+a+b+d)^2 (5+a+b+c+d)^2
(1+d)(2+d)}\\
&\quad {\scriptstyle \times (3+a+d)(4+a+c+d)(2+a)(3+a+c)(1+a)(2+a+c)(1+c)(9+2a+2b+c+2d)}\\
 &\quad {\scriptstyle \times (8+2a+2b+c+d)(7+2a+2b+c+d)(6+a+2b+c+d)(5+a+2b+d)(7+2a+b+c+d)}\\ 
&\quad {\scriptstyle \times (6+2a+b+c+d)
(5+2a+b+c)(4+a+b+c)(3+a+b)(3+a+b+c)(2+a+b)}
\\
\tdim V_{a\o_6 +b\o_1}&=
{\scriptstyle \frac{1}{548674560000}
 (1+b)(2+b)(3+b)(4+b)(5+b)(9+a+b)(8+a+b)(7+a+b)}\\
 & \quad {\scriptstyle \times
(6+a+b)(5+a+b)(7+a)(6+a)(5+a)^2(4+a)^2
(3+a)^2(2+a)(1+a)}.
\end{align*}

\rem This raises an obvious question. To what extent are the dimension formulas 
proved in \cite{LMadv,LMuniv}, valid for intermediate Lie algebras? In particular,
 in \cite{LMuniv} we gave a general dimension formula for the Cartan powers of a 
simple Lie algebra $\fg$  in terms of its Vogel's parameters $\a,\b,\g$. 
Theorem \ref{dimexc} is the specialization of that formula to the exceptional
series, and extends to the intermediate Lie algebra $\fe_{7\frac{1}{2}}$
with Vogel's parameters $\a=-2$, $\b=10$, $\g=16$. 

Also, the remark we made at the end 
of section 5 shows that the  formula
for $\tdim\fg\up k$  holds for $\widetilde{\fsl}_n$ with the same parameters 
$\a=-2$, $\b=2$, $\g=n$ as for $\fsl_n$. 
Another interesting case is the intermediate Lie algebra of $\fsp_{2n+2}$, the odd 
symplectic algebra $$\fsp_{2n+1}=\fsp_{2n}\op\CC^{2n}\op\CC.$$
We have seen that as an $\fsp_{2n}$-module, $\fsp_{2n+1}^{(k)}\simeq 
S^{2k}(\CC^{2n}\op\CC)$, which has dimension $\binom{2n+2k}{2k}$. 
Again, that's exactly what our dimension formula predicts for 
Vogel's parameters $\a=-2$, $\b=1$, $\g=n+\frac{5}{2}$.
 
\smallskip 
\noindent{\it Question}:
How could one incorporate the intermediate Lie algebras into the 
formalism of the universal Lie algebra developped by Vogel and Deligne? A first   
obstacle is that we not longer have an invariant quadratic form, which was  a
basic ingredient in their categorical constructions.

 \section{Sextonionic geometry}\label{geomsect}

In this section we study a few projective varieties that can be defined 
naturally in terms of the sextonions, in the same way as some more familiar 
varieties are defined in terms of the usual (complexified) 
composition algebras. In particular,we 
investigate in some detail the geometry of the projective plane over $\SS$, 
which is a singular but close cousin of the famous four Severi varieties
$\AA\PP^2$, for  $\AA = \RR,\CC,\HH, \OO$. 
Then we consider the Grassmannian $G_{\o}(\SS^3,\SS^6)$, again a singular 
variety but which shares the very nice geometric properties of the 
smooth varieties $G_{\o}(\AA^3,\AA^6)$ for  $\AA = 0, \RR,\CC,\HH, \OO$
\cite{LMmagic}.

\subsection{$\SS$-lines} For $\AA=\RR,\CC,\HH,\OO$, an $\AA$-line is a smooth quadric
of dimension $a$, and can be described as the image of the Veronese map
$$\nu_2 : \BP (\AA\op\AA )\dashrightarrow\PP\cJ_2(\AA), \qquad \nu_2(x,y)=\begin{pmatrix} x\bar{x}  &  x\bar{y} \\
y\bar{x}  &  y\bar{y}   \end{pmatrix}.$$
Here $\cJ_k(\AA)$ denotes the algebra of Hermitian matrices of order $k$ with coefficients
in $\AA$. The image of this map is the quadric defined by the vanishing of the determinant.

All this makes perfect sense for $\AA=\SS$, except that the determinantal quadric in $\PP\cJ_2(\SS)$
is not smooth. Indeed, $\cJ_2(\SS)=\cJ_2(\HH)\op\cA_2(\BH^{\perp})$, where $\cA_2(\BH^{\perp})$ denotes
the (two-dimensional) space of skew-symmetric matrices with coefficients in $\BH^{\perp}\subset \BS$.
If we write a matrix $M\in\cJ_2(\SS)$ as $M=R+S$, with $R\in\cJ_2(\HH)$ and $S\in \cA_2(\BH^{\perp})$,
 then $\det(M)=\det(R)$. We conclude that:

\begin{quote} {\it An $\SS$-line $\SS\PP^1$ is a singular quadric of dimension $6$ in
$\PP\cJ_2(\SS)\simeq\PP^7$, singular along a line. }
\end{quote}

\subsection{The sextonionic plane} 
For $\AA=\RR,\CC,\HH,\OO$, the $\AA$-plane $\AA\PP^2$
can be defined as  the image of the Veronese map
$$\nu_2 : \BP (\AA\op\AA\op\AA )\dashrightarrow\PP\cJ_3(\AA), \qquad \nu_2(x,y,z)=
\begin{pmatrix} x\bar{x}  &  x\bar{y} &  x\bar{z} \\
y\bar{x}  &  y\bar{y}  &  y\bar{z} \\
z\bar{x}  &  z\bar{y}  &  z\bar{z} \end{pmatrix}.$$

While  $\cJ_3(\SS)$ is a Jordan
algebra, in fact a Jordan subalgebra of the exceptional simple Jordan algebra
$\cJ_3(\OO)$, it is {\it not simple}. In fact, we can write
$\cJ_3(\SS)=\cJ_3(\HH)\op\cA_3(\BH^{\perp})$. A   computation shows that
$\cA_3(\BH^{\perp})$ is a two-sided Jordan ideal of $\cJ_3(\SS)$, and its square is
obviously zero. Therefore, $\cA_3(\BH^{\perp})$ is the radical of $\cJ_3(\SS)$,
whose semi-simple part is $\cJ_3(\HH)$.

\begin{proposition}
The derivation algebra of $\cJ_3(\SS)$ is $Der\cJ_3(\SS)\simeq\fg(\RR,\SS)$.
\end{proposition}

The same statement holds for the normed algebras, and the proof of \cite{bs}
works for $\SS$ without change. \medskip

Let $x,y,z\in\HH$ and $r,s,t\in\BH^{\perp}$. Then
$$\nu_2(x+r,y+s,z+t)=
\begin{pmatrix} x\bar{x}  &  x\bar{y} &  x\bar{z} \\
y\bar{x}  &  y\bar{y}  &  y\bar{z} \\
z\bar{x}  &  z\bar{y}  &  z\bar{z} \end{pmatrix} +
\begin{pmatrix} 0  &  r\bar{y}-xs &  r\bar{z}-xt \\
s\bar{x}-yr  &  0  &  s\bar{z}-yt \\
t\bar{x}-zr  & t\bar{y}-zs  &  0 \end{pmatrix}.$$
The first summand is in $\cJ_3(\HH)$, and the second in $\cA_3(\BH^{\perp})$,
since $\BH^{\perp}$ is a two-sided ideal of $\SS$.
\smallskip

Now recall that there is a natural identification of $\cJ_3(\HH)$ with $\La 2\CC^6$,
such that the $\HH$-plane $\HH\PP^2$ is identified with the Grassmannian $G(2,6)\subset
\PP\La 2\CC^6$. Let $W=\CC^6$.

\begin{proposition}\label{s2inc}
There is a natural identification of $\cJ_3(\SS)$ with $V=\La 2W\op W^*$, such
that $\SS\PP^2$ is identified with the closure of the set
of pairs $[\sigma,w]\in\PP V$, where $\sigma$ belongs to $G(2,6)$,
and $w$ represents a hyperplane of $W$ containing the plane $\sigma$.
\end{proposition}

\proof The fact that $\cJ_3(\SS)$ may be identified  with $V=\La 2W\op W^*$ 
was noticed in 6.2. Now $\SS\PP^2$ is a subvariety of $\PP V$, stable under
the natural action of the intermediate Lie algebra $\fg=\fsl_6\op\La 3\CC^6\op\CC$. 
An easy explicit computation shows that it contains the set of pairs
$[\sigma,w]$, where $\s$ represents a plane contained in the hyperplane defined
by $w$. But this is a rank-four vector bundle over $G(2,6)$, hence an irreducible 
variety of dimension $12$, hence an open subset of the irreducible variety $\SS\PP^2$. 
This implies our claim. \qed 

\begin{corollary}
The variety  $\SS\pp 2$  is singular along $\PP W^*\simeq\PP^5$.
\end{corollary}

This is in   agreement with the principle stated in \cite{chaput2}, following
which the very nice algebraic properties of the normed algebras have their geometric
counterpart in the smoothness of the associated projective varieties. For example,
$\cJ_k(\OO)$ is no longer a Jordan algebra for $k\ge 4$, and every natural definition
of the $\OO$-projective space $\OO\PP^{k-1}$ gives a singular variety.

\begin{corollary}
The action of $PSL_6$ on $\SS\PP^2$ has three orbits:  the singular
locus $\PP W^*$, the Grassmannian $G(2,6)\subset\PP \La 2W$, and their complement.
The smooth locus of $\SS\PP^2$ is the total space of a rank four homogeneous vector
bundle over $G(2,6)$.
\end{corollary}

\medskip
The projective planes $\AA\PP^2\subset\PP\cJ_3(\AA)$ are the four Severi varieties,
the only smooth $n$-dimensional varieties $X\subset\PP^m$, with $m=\frac{3n}{2}+2$, whose secant
variety (the determinantal cubic) is not the whole ambient space. The $\SS$-plane has the same
properties, except that it is not smooth, as we have just seen. (Note that,
$\BS\pp 2$ is not optimal for Zak's theorem on singular varieties with
secant defect, see \cite{zak}, II.2.8, although it is naturally contained 
in $J(\BP W^*,\BH\pp 2)$
defined below, which is optimal, and the two varieties have the same secant variety.)

\begin{proposition}
The secant variety of $\SS\PP^2$ is the determinantal cubic, a cone over the
determinantal cubic in $\PP\cJ_3(\HH)$.
\end{proposition}

\proof The secant variety is clearly contained in the determinantal cubic. Equality
means that any pair $(\o,h)$, where $\o\in\La 2W$ has rank four and $h$ is a generic
linear form, can be written as a sum $(\a,k)+(\b,l)$, where $\a,\b$ have rank two,
and $k$ (respectively $l$) defines a hyperplane containing the plane $A$ (respectively
$B$) defined by $\a$ (respectively $\b$). This implies that $k_{|B}=h_{|B}$ and
$l_{|A}=h_{|A}$. Conversely, we can choose any decomposition $\o=\a+\b$ into a sum
of rank two elements, define $k$ and $l$ on $A\op B$ by the conditions that
$k_{|A}=0$, $k_{|B}=h_{|B}$ and $l_{|A}=h_{|A}$, $l_{|B}=0$, and then adjust
freely on a complement $C$ of $A\op B$ so that $(k+l)_{|C}=h_{|C}$. Then $h=k+l$,
and we are done. \qed

\subsection{Orbits in $\cJ_3(\BS)$}

\begin{proposition}
The action of $PSL_6$ on $\PP V$ is prehomogeneous.
The open orbit is the complement of the determinantal cubic and
the linear subspace $\PP\La 2W\simeq\PP^{14}$.
In fact there are exactly nine $PSL_6$-orbits in $\PP V$.
\end{proposition}

The orbits are very easy to describe. For a pair $(\o,h)\in\La 2W$, the rank of
$\o$ can be $0$, $2$, $4$ or $6$, and $h$ can define a hyperplane containing or
not the kernel of $\o$, or be zero.
The incidence diagram is as follows, where $\cO_k$ denotes an
orbit of dimension $k$:
$$\begin{array}{ccccccccc}
 & & \cO_{19} & \ra & \cO_{14} & & & &  \\
 &  \nearrow  & & \searrow  & & \searrow & & & \\
\cO_{20} & & &  & \cO_{17} & \ra & \cO_{12} & \ra & \cO_{5}  \\
 &  \searrow  & & & & \searrow  & & \searrow &  \\
& & \cO'_{14} & \ra & \ra & \ra & \cO_{13}& \ra & \cO_8
\end{array}$$

It is more natural to consider the action on $\PP\cJ_3(\SS)$
of the automorphism group of $\cJ_3(\SS)$, or of the group
$PSL(3,\SS)$ preserving the determinant. We define $SL(3,\SS)$ to be the closed
subgroup of $GL(\cJ_3(\SS))$ with Lie algebra $\fg(\CC,\SS)\simeq Der(\cJ_3(\SS))\op
\cJ_3(\SS)_0$, where the space of traceless matrices $\cJ_3(\SS)_0$ acts on $\cJ_3(\SS)$
by multiplication. (Recall that the Lie algebra structure follows from the fact that for
any  $x,y\in\cJ_3(\SS)$, the bracket $D_{x,y}=[M_x,M_y]$ of the multiplication operators
by $x$ and by $y$, is a derivation of $\cJ_3(\SS)$.) In fact, $\fg(\CC,\SS)$ is our intermediate
Lie algebra $\fg$. 

Clearly, $G(2,6)$ is {\it not} stable under the action of $PSL(3,\SS)$, since otherwise
$\La 2\CC^6$ would be stable under the action of $\fg(\CC,\SS)$. We conclude:

\begin{proposition}
The action of $PSL(3,\SS)$ on $\SS\PP^2$ has only two orbits:  the singular
locus $\PP W^*$, and the smooth locus.
\end{proposition}

We now examine the $PSL(3,\BS)$ orbits in $\BP \cJ_3(\BS)$. Let
$J(\BP W^*, \BH\BP^2)$ denote the cone over $\HH\PP^2=G(2,6)$, and note
that $\BS\pp 2\subset J(\BP W^*, \BH\BP^2)$.
\def\coneh{J(\BP W^*, \BH\BP^2)}
A point $p\in J(\BP W^*, \BH\BP^2)\backslash \SS\PP^2$,  can be represented by a sum
$$ \begin{pmatrix} 1  &  0 &  0 \\
 0 & 0  & 0 \\ 0  & 0  &  0 \end{pmatrix}+\begin{pmatrix} 0  &  r &  s \\
-r  &  0  &  t \\  -s& -t  &  0 \end{pmatrix}$$
with $t\ne 0$. We prove that the tangent space $T_p=\fg(\SS,\CC).p$ to the orbit of this point,
has the same dimension as the cone over $\HH\PP^2$, implying that
$\coneh\backslash\BS\pp 2$ is a single $PSL(3,\BS)$ orbit.

First note that the action of $PSL(3,\HH)$ contributes by the dimension of $\HH\PP^2$.
What remains to prove is that $\cA_3(\BH^{\perp})$ is contained in $T_p$. Recall that
$\fg(\SS,\CC)=Der\cJ_3(\SS)\op \cJ_3(\SS)_0$. By the left action of $\cA_3(\BH^{\perp})
\subset\cJ_3(\SS)_0$, we get that
$$\begin{pmatrix} 0  &  a &  b \\
-a  &  0  &  0 \\  -b& 0  &  0 \end{pmatrix}\in T_p \qquad\forall a,b\in\BH^{\perp}.$$
Then we use the action of the triality algebra $\ft(\HH)\subset\ft(\SS)\subset Der\cJ_3(\SS)$.
Recall that $\ft(\HH)\simeq\fsl(A)\times\fsl(B)\times\fsl(C)$, where $A,B,C$ have dimension
two, and that $\SS\simeq A\ot B\op C$. In particular, $t$ being a nonzero vector in $C\simeq\BH^{\perp}$ can be taken to any vector in $C$, so that
$$\begin{pmatrix} 0  &  0 &  0 \\
0  &  0  &  c \\  0 & -c  &  0 \end{pmatrix}\in T_p \qquad\forall c\in\BH^{\perp}.$$
Thus $\cA_3(\BH^{\perp})\subset T_p$, and our claim is proved.

\medskip
Now, a point in   $\s(\BS\pp 2)\backslash \coneh$,  can be represented by a sum
$$ \begin{pmatrix} 1  &  0 &  0 \\
 0 & 1  & 0 \\ 0  & 0  &  0 \end{pmatrix}+\begin{pmatrix} 0  &  r &  s \\
-r  &  0  &  t \\  -s& -t  &  0 \end{pmatrix}.$$
The action of $\cA_3(\BH^{\perp})$ by multiplication is trivial on the second factor.
Since for $a,b,c\in\BH^{\perp}$,
$$\begin{pmatrix} 0  &  a &  b \\
-a  &  0  &  c \\  -b& -c  &  0 \end{pmatrix} \begin{pmatrix} 1  &  0 &  0 \\
 0 & 1  & 0 \\ 0  & 0  &  0 \end{pmatrix}=\begin{pmatrix} 0  &  a &  b/2 \\
-a  &  0  &  c/2 \\  -b/2& -c/2  &  0 \end{pmatrix},$$
we see that the orbit of this point must be open in the determinantal hypersurface
$\s (\BS\pp 2)$,
independantly of $r,s,t$. We conclude:

\begin{proposition}
The orbit closures of $PSL(3,\SS)$ in $\PP\cJ_3(\SS)$ are
$$\PP^5=\SS\PP^2_{sing}\subset \SS\PP^2\subset\coneh
\subset \s(\BS\pp 2) \subset\PP\cJ_3(\SS).$$
 \end{proposition}

As in the case of  $\AA\PP^2\subset\PP\cJ_3(\AA)$, with $\AA$ normed, we get
a simple chain of orbit closures. The two differences here are
that $\SS\PP^2$ is singular, and a proper subvariety of the cone over $\HH\PP^2$.

\subsection{Linear spaces in $\SS\PP^2$}

\begin{proposition}
The open orbit of the adjoint variety $X^{ad}(\SS,\CC)$ parametrizes a family of 
$\PP^4$'s in $\SS\PP^2$. This family has dimension $15$.
\end{proposition}

\proof We prove that if $x$ belongs to the open orbit in $X^{ad}(\SS,\CC)$,
its image in $\cJ_3(\SS)$ defines a $\PP^4$ contained in $\SS\PP^2$.
By homogeneity, we may suppose that $x$ belongs to the adjoint variety $\FF_{1,5}$ 
of $\fsl_6$, and corresponds to a pair $(\ell\subset H)$, for $\ell$ a line and $H$ 
a hyperplane in $W\simeq\CC^6$. Its action on $\cJ_3(\SS)=\La 2W\op W^*$ has for image 
$\ell\we H\op H$, a five dimensional vector space. The projectivization of this vector
space is clearly contained in $\SS\PP^2$, by Proposition \ref{s2inc}, 
because a nonzero vector in $\ell\we H$ defines a two-plane containing $\ell$ 
and contained in $H$. \qed 

%\medskip A rather unexpected fact is that we have other families of $\PP^4$'s
%on $\SS\PP^2$. In fact it is easy to describe all the linear spaces on this variety. 
%Of course we just need to understand the unextendable ones. 

\medskip Unlike the Severi varieties, there are other families of unextendable
linear spaces on $\BS\pp 2$:
 
\begin{proposition}
The unextendable linear spaces on $\BS\pp 2$ are as follows:
\begin{itemize}
\item A $10$-dimensional family of $\BP^5$'s parametrized by a one
point compactification of a line bundle over $G(3,W)$.
\item An irreducible family of $\BP^4$'s of dimension $15$,
with an open subset  
given   by the smooth locus of the adjoint variety $X^{ad}(\SS,\CC)$, 
 
 \item An irreducible family of $\BP^4$'s of dimension $15$,
with an open subset  
given   by  
  the total space of the vector bundle $\La 2Q^*(1)$, where $Q$ denotes the rank 
four quotient bundle on $\PP W$. 
\end{itemize}
\end{proposition}

\def\oP{\overline P}

\proof Let $P\subset\SS\PP^2$ be an unextendable linear space. Its projection to $\PP\La 2W$
is a linear space contained in $G(2,6)$, so is either the 
set of planes containing a line $\ell$ and contained in a $k$-dimensional
space $L$, or the set of planes contained in a three plane $M$. 

In the second case, again in an adapted basis, $P$ must be generated by vectors
$e_1\we e_2+ze_3^*,  e_2\we e_3+ze_1^*,  e_3\we e_1+ze_2^*, e_4^*,  e_5^*,  e_6^*$.
We thus get a family of $\PP^5$'s on $\SS\PP^2$, parametrized by a $\CC$-bundle 
over the Grassmannian $G(3,6)$. This family becomes complete when we add to it a single 
point, corresponding to the singular set $\PP W\simeq\PP^5$ of $\SS\PP^2$. 

In the first case, in an adapted basis, $P$ must be generated by vectors 
$e_1\we e_2+h_2,\ldots, e_1\we e_{k}+h_{k}, e_{k+1}^*,\ldots ,e_6^*$, where $h_2,
\ldots ,h_{k+2}$ are linear forms such that $h_i(e_1)=0$ for all $i$ and 
the matrix $h_i(e_j)$, $2\le i,j\le k+2$, is skew-symmetric. In particular, 
$P$ has affine dimension $5$. If $k=0$, we get the singular $\PP^5$. Note also that $k\ne 1,2$, 
otherwise $P$ would be extendable. Thus $k\ge 3$ and we get a family 
of $\PP^4$'s on $\SS\PP^2$ of dimension $k(6-k)+(k-1)+\binom{k-1}{2}=k(6-k)+\binom{k}{2}$
(choice of $L$ plus choice of $\ell\subset L$ plus choice of $h$). 
Note that for $k=5$, we recover the $15$-dimensional family parametrized by the open
orbit of the adjoint variety. But $k=6$ gives another family of the same dimension. 
Note that in that case, $\ell$ being the line generated by $e_1$, the map $h$ should
be seen as a skew-symmetric morphism from $W/\ell$ to $\ell^{\perp}\simeq (W/\ell)^*$,
depending linearly on the vector we choose on $\ell$. Thus 
our family is parametrized by the vector bundle $\La 2Q(1)$ on $\PP W$. 

Finally, it is easy to check that the 
other cases belong to the closure of these two maximal families. \qed

\subsection{Point-line geometry}
When $\AA$ is a normed algebra, the $\AA$-plane is covered by a family of $\AA$-lines
(i.e. $\BA\pp 1$'s) parametrized by
$\AA\PP^2$ itself. This family of $\BA\pp 1$'s defines a plane projective geometry
on $\AA\PP^2$, in the sense that two generic points are joined by a unique line,
and two generic lines meet in a unique point.
We now show that the same picture   holds for $\SS\PP^2$. 

The $\AA$-lines can be described as the entry-loci of the points inside the secant cubic.
For the sextonions, we choose a pair $(\o,h)$, where $\o\in\La 2W$ has rank four, and $h$
is a linear form. Denote by  $P$ the support of $\o$, i.e., the four plane
which is the image of the contraction by $W^*$. A computation shows that the entry-locus
of $(\o,h)$ is the intersection of $\SS\PP^2$ with the linear space $\phi_h(\La 2P)+ P^{\perp}$,
where $\phi_h=Id_{\La 2W}+\psi_h$ for an endomorphism  $\psi_h : \La 2P\ra P^*$ defined
by $h$, more precisely by the   restriction of $h$ to $P$.

To be more explicit, suppose that $\o=e_1\we e_2+e_3\we e_4$, and let's try to
solve the equation $(\o,h)=(\a,k)+(\b,l)$. Around $\a_0=e_1\we e_2$, $\b_0=e_3\we e_4$,
a solution of the equation $\o=\a+\b$ can be written
$$\begin{array}{rcl}
\a & =& (1+sv-ut)^{-1}(e_1+se_3+te_4)\we (e_2+ue_3+ve_4),\\
\b & = & (1+sv-ut)^{-1}(e_3-ve_1+te_2)\we (e_4+ue_1-se_2).
\end{array}$$
Then $h=k+l$, with $k_{|\a}=0$ and $l_{|\b}=0$, if $k(e_1+se_3+te_4)=k(e_2+ue_3+ve_4)=0$,
$k(e_3-ve_1+te_2)=h(e_3-ve_1+te_2)$ and $k(e_4+ue_1-se_2)=h(e_4+ue_1-se_2)$. This gives
$$\begin{array}{rcl}
(1+sv-ut)k(e_1) & = & (sv-ut)h(e_1)-sh(e_3)-th(e_4), \\
(1+sv-ut)k(e_2) & = & (sv-ut)h(e_2)-uh(e_3)-vh(e_4), \\
(1+sv-ut)k(e_3) & = & h(e_3)-vh(e_1)+th(e_2), \\
(1+sv-ut)k(e_4) & = & h(e_4)+uh(e_1)-sh(e_2).
\end{array}$$
Letting $h(e_i)=h_i$, and completing $e_1, e_2, e_3, e_4$ into a basis of $W$,
we deduce that the projective span of $(\a,k)$ is
$$\begin{array}{l}
[\sum_{1\le i<j\le 4}Z_{ij}e_i\we e_j,\;
h_1Z_{34}+h_3Z_{23}+h_4Z_{24},\;
h_2Z_{34}+h_3Z_{13}-h_4Z_{14}, \\
\hspace{3.45cm}
h_3Z_{12}-h_1Z_{14}-h_2Z_{24},\;
h_4Z_{12}+h_1Z_{13}+h_2Z_{23},\; k_5,k_6],
\end{array}$$
and this is a point of $\SS\PP^2$ when $\sum_{1\le i<j\le 4}Z_{ij}e_i\we e_j$ has rank two.
We can describe this linear space in a more invariant way as follows. The two form $\o$
defines the four-space $P$, and the non-zero vector $\we^2\o\in\La 4P$ which allows 
one to
identify $P$ with its dual. Choose a supplement $P^{\circ}$ to $P^{\perp}$ in $W^*$, so
that the composition $P^{\circ}\hookrightarrow W^*\ra P^*$ is a natural isomorphism.
Then the linear space above is the set of vectors
$$Z + Z\we (h\lrcorner\o)/(\o\we\o)+ Y, \qquad Z\in\La 2P, \quad Y\in P^{\perp},$$
where $h\lrcorner\o$ belongs to $P$, hence $Z\we (h\lrcorner\o)$ to $\La 3P=P^*\otimes\La 4P$,
so that we obtain after division by $(\o\we\o)$ a vector in $P^*$ that we identify
with $P^{\circ}\subset W^*$. The resulting vector is uniquely defined only up to $P^{\perp}$,
but the $Y$ term allows one to ignore that point.

This space is therefore defined only by the tensor $\o\we\o + h\lrcorner\o\in\La 4W\op W$,
where $\o\we\o$ is a decomposable tensor in $\La 4W\simeq \La 2W^*$, and $h\lrcorner\o\in W$
is a linear form on $W^*$ vanishing on the plane defined by $\o\we\o$. We finally get a family
of $\SS$-lines parametrized by the smooth part of the dual plane $\SS {\Hat \PP}^2$.

It remains to understand how these $\SS$-lines degenerate when we approach the singular set
of this dual plane. To see this, we compute the entry locus of a generic point of the
determinantal hypersurface in $\PP\cJ_3(\SS)$, of the form $\o+h$, with $\o\in\La 2W$ a
decomposable tensor and $h$  linear form which is not identically zero on the plane defined by
$\o$. We check that this entry locus only depends on the kernel of the restriction of $h$
to that plane: precisely, if $e$ is a generator of that line, it is a smooth $8$-dimensional
quadric obtained as the intersection
of $\SS\PP^2$ with the linear space $e\we W\op e^{\perp}$. Such a smooth quadric
is clearly covered by $6$-dimensional quadrics singular along a line, which can be
obtained as limits of $\SS\PP^1$'s on $\SS\PP^2$. And the family of these smooth
quadrics is naturally parametrized by $\PP W$, the singular set of the dual plane
$\SS {\Hat \PP}^2$.

\medskip

Using this explicit description, we easily get:

\begin{proposition}
Two generic $\SS$-lines on $\SS\PP^2$ meet in a unique point. Through two
generic points of $\SS\PP^2$ passes a unique $\SS$-line.
\end{proposition}

\subsection{The first row}
To pass to the variety $X= \BS\BP^2_0\subset \BP V$ of the first row, as with
the rest of the series we take a hyperplane section, but now the
hyperplane section is no longer generic, as it cuts only the
first factor. The variety $\BS\BP^2_0$ has a corresponding
description where $G(2,W)$ is replace by
the $\o$-isotropic Grassmanian $G_{\o}(2,W)$.

\subsection{The  Grassmannian $G_{\o}(\SS^3,\SS^6)$}
The varieties from the third line of the geometric Freudenthal square have several
interesting interpretations, as Lagrangian Grassmannians of symplectic $\AA$-subspaces of $\AA^6$,
or cubic curves over the simple Jordan algebras $\cJ_3(\AA)$, or conformal compactifications
of these Jordan algebras. Recall that they are defined as the closures of the images of the
maps
$$\nu_3 : \cJ_3(\AA)\rightarrow\PP\cZ_2(\AA), \qquad \nu_3(x)=
\begin{pmatrix} 1  &  x    \\
Q(x)  &  \det(x)  \end{pmatrix},$$
where $Q(x)$ denotes the cofactor matrix of $x$ (see \cite{LMmagic},
section 1.2, or \cite{clerc}, section 6). We use the same definition over the
sextonions. Our first claim is about the equations of the resulting variety
$G_{\o}(\SS^3,\SS^6)$. (For the nondegenerate case, this is Proposition
6.2 in \cite{clerc}, but the proof is not quite correct). The following argument works
in general. We begin by exhibiting a set of quadratic equations of $G_{\o}(\SS^3,\SS^6)$, 
which define it set-theoretically.
% (but beware that we  do not get the whole set of quadratic 
%equations). 
 
\begin{lemma}\label{gseqn}
The variety $G_{\o}(\SS^3,\SS^6)\subset \PP\cZ_2(\SS)$  is the set of matrices
$\begin{pmatrix} s  &  x  \\ y  & t  \end{pmatrix}$, such that
$$Q(x)=sy, \quad Q(y)=tx, \quad xy=st I.$$
\end{lemma}

\proof We must prove that such a matrix belongs to the closure of $\nu_3(\cJ_3(\SS))$.
This is clear if $s\ne 0$. Since
$$\nu_3(x^{-1})=\begin{pmatrix} \det(x)  &  Q(x) \\ x  & 1  \end{pmatrix},$$
this is also true for $t\ne 0$.
But $w\in \cJ_3(\SS)$ acts on $\cZ_2(\SS)$ by the translation
$$t_w\begin{pmatrix} s  &  x  \\ y  & t  \end{pmatrix}
=\begin{pmatrix} s  &  x+sw  \\ y+2Q(x,w)+sQ(w)  & t+\trace(yw)+\trace(xQ(w))+s\det(w)
\end{pmatrix},$$
where $Q(x,w)$ denotes the polarization of $Q$.
This action of $\cJ_3(\SS)$
preserves our set of quadratic equations, but clearly not the subspace
of matrices such that $t=0$. The claim follows. \qed

\medskip 
Let $Sp(6,\BS)$ denote the closed subgroup of $GL(\cZ_2(\SS))$ 
 defined by the Lie algebra $\fg(\SS,\HH)$. 
It contains the group $Sp(6,\HH)=Spin_{12}$, whose action on
$\cZ_2(\SS)$ leaves  invariant the subspace
$\cZ_2(\HH)\simeq \Delta_+$. Remember that the closed orbit of 
$PSp(6,\HH)$ in $\PP\cZ_2(\HH)$ is the spinor variety 
$\SS_+=G_{\o}(\HH^3,\HH^6)$.

\begin{proposition}
The variety $G_{\o}(\SS^3,\SS^6)\subset\PP\cZ_2(\SS)$ is the highest weight
variety of $PSp(6,\BS)$ in $\BP\cZ_2(\SS)$.
It can be interpreted as the subvariety of
$\PP(\Delta_+\op U)$,  consisting of the closure of the set of pairs 
$(\Sigma, u)$ where
$\Sigma$ defines a maximal isotropic subspace in the family 
$\SS_+=G_{\o}(\HH^3,\HH^6)$, and $u$ belongs to $\Sigma$.
\end{proposition}

\proof 
A stable complement is given by the space of matrices
with coefficients in $\BH^{\perp}$. This complement has dimension 
$12$ and, since the
action is nontrivial, it must coincide with the natural 
representation of $Spin_{12}$
on $U\simeq\CC^{12}$. (This proves the description of $V$ in 
\ref{thirdrowsubsect}.)

We have $S^2(\Delta_+\op U)^*=S^2\Delta_+\op (\Delta_+\ot U)\op S^2U$. Using the notation
of Bourbaki for the weights of $\fso_{12}$, we have $\Delta_+=V_{\o_6}$ and
$$S^2\Delta_+=V_{2\o_6}\op V_{\o_2}, \quad  \Delta_+\ot U=V_{\o_1+\o_6}\op V_{\o_5},
\quad S^2U=V_{2\o_1}\op\CC.$$
Now recall Lemma \ref{gseqn}, and decompose $\cZ_2(\SS)$ into $\cZ_2(\HH)\op (\SS^{\perp}
\op\SS^{\perp})=\Delta_+\op U$. We see that $G_{\o}(\SS^3,\SS^6)$ has
quadratic equations of different types: those only involving $\cZ_2(\HH)$ are the 
quadratic equations of $G_{\o}(\HH^3,\HH^6)$, which gives $V_{\o_2}$. There are also
equations of mixed type,  i.e. from $\Delta_+\ot U$. 

It follows that $G_{\o}(\SS^3,\SS^6)$ can be
defined as a set pairs $(\Sigma,u)$, where $\Sigma\in\PP\Delta_+$ defines a maximal
isotropic subspace of $U$ in one of the two families of these, and $u$ belongs to
some subspace of $U$ defined by $\Sigma$ in some invariant way. The only possibility
is that this space is $\Sigma$ itself (it cannot be zero since $G_{\o}(\SS^3,\SS^6)$ 
is certainly not contained in $G_{\o}(\HH^3,\HH^6)$, and it cannot be the whole of $U$ 
since  we do have mixed equations). In particular, $u$ must be isotropic, and the
trivial factor of $S^2U$ must appear in the space of quadratic equations of
$G_{\o}(\SS^3,\SS^6)$. 

\rem We conclude that the space of quadratic equations of $G_{\o}(\SS^3,\SS^6)$, as 
an $\fso_{12}$-module, is $V_{\o_2}\op V_{\o_5}\op \BC=\fso_{12}\op\Delta_-\op\CC$. 
But this is just the 
intermediate Lie algebra $\fg=\fg(\SS,\HH)$, which is no surprise since on the 
third line of the magic chart, we have the invariant symplectic form $\omega$,
which allows to associate to every vector $x\in\fg$  the quadratic form 
$q_x(v)=\omega(v,xv)$.

\begin{corollary}
$G_{\o}(\SS^3,\SS^6)$ is singular along the quadric $\QQ^{10}\subset\PP U$.
Its smooth locus has two orbits under the action of $PSO_{12}$, but is homogeneous
under the action of $PSp(6,\SS)$.
\end{corollary}

\proof The Zariski tangent space of $G_{\o}(\SS^3,\SS^6)$ at a point $(0,u)\in
\QQ^{10}\subset\PP U\subset\PP(\Delta_+\op U)$ certainly contains the line of
$\Delta_+$ generated by $\Sigma\in\PP\Delta_+$ parametrizing any maximal isotropic
subspace of $U$ containing $u$. The linear span of such $\Sigma$'s is isomorphic
with a half-spin representation of $Spin(u^{\perp}/\CC u)=Spin_{10}$ - in particular,
its dimension is $16$. But our Zariski tangent space also contains $U$, obviously,
and we already get $16+12=28$ dimensions, which is more than the dimension, $21$,
of $G_{\o}(\SS^3,\SS^6)$. This proves the first claim.

The complement of $\QQ^{10}$ is the set of pairs $(\Sigma,u)$, where $\Sigma$ is non
zero and parametrizes a maximal isotropic subspace of $U$ containing $u$.  The action
of $PSO_{12}$ gives two orbits, one where $u=0$ and one where
$u\neq 0$. But the condition $u=0$
is not $\fg(\SS,\HH)$-invariant, so the complement of $\QQ^{10}$ is
$PSp(6,\SS)$-homogeneous, hence smooth, and exactly equal to the smooth locus
of $G_{\o}(\SS^3,\SS^6)$. \qed

%The tangent variety of $G_{\o}(\SS^3,\SS^6)$ is a quartic hypersurface, whose
%complement is $PSp(6,\SS)$-homogeneous.

\begin{proposition}
The orbit closures of $PSp(6,\SS)$ in $\PP\cZ_2(\SS)$ are the cones over the four
$PSp(6,\HH)$-orbits in $\PP\cZ_2(\HH)$ with vertex $\PP U$, the variety
$G_{\o}(\SS^3,\SS^6)$ and its singular locus $\QQ^{10}\subset\PP U$.
\end{proposition}

\proof Consider a point in $\PP\cZ_2(\HH)\subset\PP\cZ_2(\SS)$, given by some matrix
$m=\begin{pmatrix} s  &  x  \\ y  & t  \end{pmatrix}$. We want to understand when the
tangent space $\fg(\SS,\HH).m$ to the orbit of $m$ contains $U\simeq\cZ_2(\BH^{\perp})$,
the subspace of $\cZ_2(\SS)$ consisting of matrices all of whose coefficients are in
$\BH^{\perp}$ (in particular, the diagonal coefficients must be  zero). Note that this
condition is certainly $PSp(6,\HH)$-invariant.

The action of $\cA_3(\BH^{\perp})$ and its dual provide us with the matrices
$$\begin{pmatrix} 0  &  su  \\ xu  & 0  \end{pmatrix}\quad {\rm and}\quad
\begin{pmatrix} 0  &  yv  \\ tv  & 0  \end{pmatrix}, \quad u,v\in\cA_3(\BH^{\perp}).$$
We can certainly solve the equations $su+yv=p$, $xu+tv=q$, as soon as the matrix
$stI-xy$ is invertible. This is the case if $m=\begin{pmatrix} 1 & 0 \\ 0 & 1 \end{pmatrix}$,
a point in the complement of the tangent quartic, which is an open $PSp(6,\HH)$-orbit in
$\PP\cZ_2(\HH)$.

This is no longer true if we consider the point $m=\begin{pmatrix} 1 & x \\ 0 & 0 \end{pmatrix}$
on the tangent hypersurface. But let us move this point by the translation $t_w$, to get the
point $$t_w(m)=
\begin{pmatrix} 1  &  x+w  \\ Q(x,w)+Q(w)  & \trace(xQ(w))+\det(w)  \end{pmatrix}.$$
Now the matrix we want to be invertible is $z_w=\trace(xQ(w))I-xQ(w)-xQ(x,w)-wQ(x,w)$.
It is enough to find $w$ such that the degree two part $wQ(x,w)$ is invertible. We claim
that this is possible as soon as the rank of $x$ is at least two. Indeed, if we represent
$x$ by some diagonal matrix with at least two nonzero eigenvalues, and if we also choose
$w$ to be diagonal, a straightforward computation shows that $wQ(x,w)$ is again diagonal
with generically nonzero eigenvalues.

We conclude that for any point in $\PP(\Delta_+\op U)$ of the form $p=(\Sigma, u)$, where
$\Sigma$ does not belong to the cone over $\SS_+=G_{\o}(\HH^3,\HH^6)$, the
$PSp(6,\SS)$-orbit of $p$ must be the whole cone over the $PSp(6,\HH)$-orbit of $\Sigma$,
with vertex $\PP U$.

A similar computation shows that when $\Sigma$ is a nonzero vector in the cone over $\SS_+$,
there are only two cases up to the $PSp(6,\SS)$-action: either $u$ does belong to $\Sigma$,
or not. \qed

\begin{proposition}
A point in $\PP\cZ_2(\SS)$, outside the tangential quartic, belongs to a unique
secant to $G_{\o}(\SS^3,\SS^6)$. In particular, $G_{\o}(\SS^3,\SS^6)$ has
only one apparent double point.
\end{proposition}

\proof We use the same argument as in \cite{clerc}: A general point
$m=\begin{pmatrix} s  &  x  \\ y  & t  \end{pmatrix}$ in $\cZ_2(\SS)$ does not
belong to the hyperplane at infinity $(s=0)$. By translation we can then suppose that $x=0$.
It is easy to check that if $y$ is invertible,
$m$ cannot belong to a secant line joining two points of $G_{\o}(\SS^3,\SS^6)$, one of
which in the hyperplane at infinity. So we need to solve the equation
$$\begin{pmatrix} s  &  x  \\ y  & t  \end{pmatrix}=
\lambda \begin{pmatrix} 1  &  a  \\ Q(a)  & \det(a)  \end{pmatrix}+
\mu \begin{pmatrix} 1  &  b  \\ Q(b)  & \det(b)  \end{pmatrix}.$$
Using the identities $Q(Q(a))=\det(a)a$ and $a Q(a)=\det(a)I$, which are valid in $\cJ_3(\SS)$,
one checks that this equation has for unique solution
$$a=\frac{\mu-\lambda}{\lambda}\frac{Q(y)}{st} \quad {\rm and}\quad
b=\frac{\lambda-\mu}{\mu}\frac{Q(y)}{st},$$
where the scalars $\lambda$ and $\mu$ are uniquely defined by the conditions that
$$\lambda+\mu=s  \quad {\rm and}\quad   \frac{(\lambda-\mu)^2}{\lambda\mu}=\frac{st^2}{\det(y)}.
\qquad\qed$$

\rem The property of having only one apparent double point is equivalent to the fact that
the projection of the variety from a general tangent space is birational. For the varieties
$G_{\o}(\AA^3,\AA^6)$, this projection can be interpreted as the map $Q :\PP\cJ_3(\AA)
\dashrightarrow\PP\cJ_3(\AA)$. This is an involutive birational isomorphism because
of the identity $Q(Q(a))=\det(a)a$, and this also holds over the sextonions.

\subsection{The adjoint variety $X^{ad}(\SS,\HH)$}
We conclude with a brief sketch of study of this variety. 
Remember the identification $\fg(\SS,\HH)=\fso_{12}\op\Delta_-\op\CC$. We define 
$X^{ad}(\SS,\HH)\subset\PP\fg(\SS,\HH)$ as the closure of the space of triples $(P,\Sigma,z)$
such that: $P\in\fso_{12}$ parametrizes 
a point of the adjoint variety $X^{ad}(\HH,\HH)$,
i.e., an isotropic plane in $U=\CC^{12}$; $\Sigma$ parametrizes a maximal isotropic space
in $U$ from the family $\SS_-$, containing $P$; $z$ is any scalar.

\begin{proposition}
The variety $X^{ad}(\SS,\HH)\subset \PP\fg(\SS,\HH)$ is the   $PSp(6,\SS)$-adjoint
variety.  Its dimension is $25$.
\end{proposition}

This is in   agreement with the fact that for the third row of Freudenthal's square,
the dimension of the adjoint variety is $4a+1$ in the nondegenerate case.

\begin{proposition}
The smooth locus of $X^{ad}(\SS,\HH)$ parametrizes a family of $8$-dimensional quadrics
on $G_{\o}(\SS^3,\SS^6)$.
\end{proposition}

We can also consider the space of lines on $G_{\o}(\SS^3,\SS^6)$. The lines which are not contained
in the singular locus form a quasi-homogeneous variety, linearly nondegenerate inside $\PP
\La {\langle 2\rangle}\cZ_2(\SS)$.

 We leave to the reader the problem of showing that
points of $X^{ad}(\SS,\HH)$, lines in $G_{\o}(\SS^3,\SS^6)$, and points of $G_{\o}(\SS^3,\SS^6)$,
are the elements -- points, lines and planes respectively, of a six-dimensional symplectic geometry.
 
 We also leave to the future the problem of studying the varieties in the sextonionic and
expanded  octonionic rows. For the octonionic row, we should get four quasi-homogeneous 
varieties defining a metasymplectic geometry in the sense of Freudenthal.

\vspace{2cm}

\begin{tabular}{lll}
Joseph M. Landsberg & \hspace*{1cm} & Laurent Manivel \\
School of Mathematics, & & Institut Fourier, UMR 5582 du CNRS \\
Georgia Institute of Technology, & & Universit\'e Grenoble I, BP 74 \\
Atlanta, GA 30332-0160 & & 38402 Saint Martin d'H\`eres cedex \\
USA & & FRANCE \\
{\rm E-mail}: jml@math.gatech.edu & &
{\rm E-mail}: Laurent.Manivel@ujf-grenoble.fr
\end{tabular}

 \end{document}